\documentclass[11pt]{article}
\usepackage{amsmath,amsthm,amssymb,amscd}

\title{The branching nerve of HDA and the Kan condition}
\author{Philippe Gaucher}

\usepackage[all,2cell]{xy}
\UseAllTwocells

\newcommand{\C}{\mathcal{C}}
\newcommand{\D}{\mathcal{D}}
\newcommand{\Z}{\mathbb{Z}}
\newcommand{\N}{\mathbb{N}}
\newcommand{\de}{\partial}
\newcommand{\p}\times

\renewcommand{\P}{\mathbb{P}}

\newcommand{\be}{\begin{equation}}
\newcommand{\ee}{\end{equation}}
\newcommand{\bea}{\begin{eqnarray}}
\newcommand{\eea}{\end{eqnarray}}
\newcommand{\beas}{\begin{eqnarray*}}
\newcommand{\eeas}{\end{eqnarray*}}

\newtheorem{thm}{Theorem}[section]
\newtheorem{prop}[thm]{Proposition}
\newtheorem{conj}[thm]{Conjecture}
\newtheorem{cor}[thm]{Corollary}
\newtheorem{nota}[thm]{Notation}
\newtheorem{defn}[thm]{Definition}
\newcommand{\bd}{\begin{defn}}
\newcommand{\ed}{\end{defn}}
\newcommand{\bp}{\begin{prop}}
\newcommand{\ep}{\end{prop}}
\newcommand{\bth}{\begin{thm}}
\renewcommand{\eth}{\end{thm}}
\newcommand{\bi}{\begin{enumerate}}
\newcommand{\ei}{\end{enumerate}}
\newcommand{\br}{\begin{rem}}
\newcommand{\er}{\end{rem}}
\newcommand{\bpf}{\begin{proof}}
\newcommand{\epf}{\end{proof}}

\newcommand{\fl}[1]{\ar@{->}[l]_{#1}}
\newcommand{\fr}[1]{\ar@{->}[r]^{#1}}
\newcommand{\fd}[1]{\ar@{->}[d]_{#1}}
\newcommand{\fu}[1]{\ar@{->}[u]^{#1}}
\newcommand{\f}[2]{\ar@{->}[#1]|{#2}}
\newcommand{\ff}[2]{\ar@2{->}[#1]|{#2}}

\newcommand{\iso}{\cong}

\newcommand{\HR}{H\!R}
\newcommand{\CR}{C\!R}
\newcommand{\HF}{H\!F}
\newcommand{\CF}{C\!F}
\newcommand{\ev}{e\!v}

\newcommand{\Rm}{\mathcal{R}^-}

\newcommand{\U}{\mathbb{U}}
\newcommand{\F}{\mathbb{F}}

\newcommand{\comp}{C\!omp(Ab)}
\newcommand{\id}{I\!d}

\date{February 2003}

\def\cartesien{%
  \ar@{-}[]+R+<6pt,-2pt>;[]+RD+<6pt,-6pt>%
  \ar@{-}[]+D+<2pt,-6pt>;[]+RD+<6pt,-6pt>%
}
\def\cocartesien{%
  \ar@{-}[]+L+<-6pt,+2pt>;[]+LU+<-6pt,+6pt>%
  \ar@{-}[]+U+<-2pt,+6pt>;[]+LU+<-6pt,+6pt>%
}

\renewcommand{\leq}{\leqslant}
\renewcommand{\geq}{\geqslant}
\renewcommand{\overline}{\widehat}

\begin{document}

\maketitle

\begin{abstract}
One can associate to any strict globular $\omega$-category three
augmented simplicial nerves called the globular nerve, the
branching and the merging semi-cubical nerves. If this strict
globular $\omega$-category is freely generated by a precubical
set, then the corresponding homology theories contain different
informations about the geometry of the higher dimensional
automaton modeled by the precubical set. Adding inverses in this
$\omega$-category to any morphism of dimension greater than $2$
and with respect to any composition laws of dimension greater
than $1$ does not change these homology theories. In such a
framework, the globular nerve always satisfies the Kan condition.
On the other hand, both branching and merging nerves never
satisfy it, except in some very particular and uninteresting
situations. In this paper, we introduce two new nerves (the
branching and merging semi-globular nerves) satisfying the Kan
condition and having conjecturally the same simplicial homology as
the branching and merging semi-cubical nerves respectively in
such framework. The latter conjecture is related to the thin
elements conjecture already introduced in our previous papers.
\end{abstract}

\tableofcontents

\section{Introduction}

An $\omega$-categorical model for higher dimensional automata (HDA)
was first proposed in \cite{Pratt}, followed by \cite{HDA} for a first
homological approach using these ideas and cubical models of
topological spaces as in \cite{Brown_cube}.

The papers \cite{Gau,sglob} demonstrate that the formalism of \textit{strict
globular $\omega$-categories} (see Definition~\ref{omega_categories})
\textit{freely generated by precubical sets} (see below)
provides a suitable framework for the introduction of new algebraic
tools devoted to the study of deformations of HDA. In particular, three
augmented simplicial nerves are introduced in our previous
papers\thinspace: the \textit{globular nerve} $\mathcal{N}^{gl}$, the
\textit{branching semi-cubical nerve} $\mathcal{N}^{-}$ and the
\textit{merging semi-cubical nerve} $\mathcal{N}^{+}$. Any
$\omega$-category freely generated by precubical sets is actually a
\textit{non-contracting $\omega$-category} (see
Definition~\ref{catnoncontract} and \cite{sglob} where a precubical
set is called semi-cubical set) and most of the theorems known so far
are expressible in this wider framework, hence the importance of the
latter notion. In this paper as well, most of the results will be
stated in the wider framework of non-contracting $\omega$-categories.

A \textit{precubical set} is a cubical set as in \cite{Brown_cube} but without
degeneracy maps of any kind.  It is easy to view such an object as a
contravariant functor from some small category $\square^{pre}$ to the
category of sets. The objects of $\square^{pre}$ are the nonnegative
integers. The small category $\square^{pre}$ is then the quotient of
the free category generated by the arrows
$\delta_i^\alpha:m\rightarrow m+1$ for $m\geq 0$, $1\leq i\leq m+1$
and $\alpha\in \{-,+\}$ by the relations $\delta_j^\beta
\delta_i^\alpha=\delta_i^\alpha \delta_{j-1}^\beta$ for $i<j$
and $\alpha\in \{-,+\}$.

Now let $K=(K_n)_{n\geq 0}$ be a precubical set. Let $I^n$ be the
$n$-dimensional cube viewed as a strict globular $\omega$-category
(see the reminder in Section~\ref{convention}). Then the
\textit{strict globular $\omega$-category freely generated by $K$} is
by definition the colimit $F(K):=\int^{n\in \square^{pre}} K_n . I^n$
where the notation $K_n.I^n$ means the coproduct of ``cardinal of $K_n$''
copies of $I^n$. This construction induces a functor from the category
of precubical sets (with the natural transformations of functors as
morphisms) to the category $\omega Cat$ of strict globular
$\omega$-categories. This functor is of course left adjoint to the
precubical nerve functor from $\omega Cat$ to the category of
precubical sets. Any strict $\omega$-category freely generated by a
precubical set is necessarily non-contracting.

Non-contracting $\omega$-categories freely generated by precubical
sets encode the algebraic properties of execution paths and of higher
dimensional homotopies between them in HDA. Indeed, the $0$-morphisms
represent the states of the HDA, the $1$-morphisms the non-constant
execution paths, and the $p$-morphisms with $p\geq 2$ the higher
dimensional homotopies between them. R. Cridlig presents
in \cite{cridlig96implementing} an
implementation with CaML of the semantics of a real concurrent
language in terms of \textit{precubical sets}, demonstrating the
relevance of this approach.

However, in such an $\omega$-category $\C=F(K)$, if $H$ is an
homotopy (that is a $2$-morphism) from an execution path (that is
a $1$-morphism) $\gamma_1$ to an execution path $\gamma_2$, then
it is natural to pose the existence of an opposite homotopy
$H^{-1}$ from $\gamma_2$ to $\gamma_1$, that is satisfying
$H^{-1} *_1 H = \gamma_2$ and $H *_1 H^{-1} = \gamma_1$. So any
$2$-morphism should be invertible with respect to $*_1$. Now let
us suppose that we are considering an $\omega$-category with a
$3$-morphism $A$ from the $2$-morphism $s_2A$ to the $2$-morphism
$t_2A$. Let $(s_2A)^{-1}$ (resp. $(t_2A)^{-1}$) be the inverse of
$s_2A$ (resp. $t_2A$) with respect to $*_1$. If there are no holes
between $s_2A$ and $t_2A$ because of $A$, then it is natural to
think that there are no holes either between $(s_2A)^{-1}$ and
$(t_2A)^{-1}$\thinspace: therefore there should exist a
$3$-morphism $A'$ such that $s_2A'=(s_2A)^{-1}$ and
$t_2A'=(t_2A)^{-1}$. Then $s_2(A *_1 A')= s_2A *_1 (s_2A)^{-1} =
s_1A$. Therefore it is natural to assume $A *_1 A'$ to be
$1$-dimensional and so the equality $A *_1 A'=s_1 A$ should hold.
This means that not only a $3$-morphism should be invertible with
respect to $*_2$, because if there is a $3$-dimensional homotopy
from $s_2A$ to $t_2A$, then there should exist a $3$-dimensional
homotopy from $t_2A$ to $s_2A$, but also with respect to $*_1$.
However the inverses of $A$ with respect to $*_1$ and $*_2$ do
not have any computer-scientific reason to be equal.

What we mean is that it is natural from a computer-scientific
point of view to deal with non-contracting $\omega$-categories
whose corresponding \textit{path $\omega$-category $\P \C$} (see
Definition~\ref{path-cat}) is a \textit{strict globular
$\omega$-groupoid} (see Definition~\ref{def_groupoid}). One would
like to emphasize the fact that it is out of the question to
suppose that the whole $\omega$-category $\C$ is an
$\omega$-groupoid because in our applications, the $1$-morphisms
of $\C$ are \textit{never} invertible, due to the irreversibility
of time.

Let $\omega Cat_1$ be the category of non-contracting
$\omega$-categories with the non-contracting $\omega$-functors as
morphisms. Let $\omega Cat_1^{Kan}$ be the full and faithful
subcategory of $\omega Cat_1$ whose objects are the non-contracting
$\omega$-categories $\C$ such that $\P\C$ is a strict globular
$\omega$-groupoid. The forgetful functor from $\omega Cat_1^{Kan}$ to
$\omega Cat_1$ commutes with limits. One can easily check that both
categories $\omega Cat_1$ and $\omega Cat_1^{Kan}$ are complete. So by
standard categorical arguments (for example the solution set condition
\cite{MR1712872}), this latter functor admits a left adjoint
$\C\mapsto \widehat{C}$ from $\omega Cat_1$ to $\omega Cat_1^{Kan}$.

What we claim is that it is then natural from a computer-scientific
point of view to deal with $\omega$-categories like $\widehat{F(K)}$
where $K$ is a precubical set. But what do the three homology theories
$H^{gl}$, $H^-$ and $H^+$ corresponding to the simplicial nerves
$\mathcal{N}^{gl}$, $\mathcal{N}^{-}$ and $\mathcal{N}^{+}$ become ?
To understand the answer given in Theorem~\ref{reponse}, we need to
recall in an informal way the thin elements conjecture which already
showed up in
\cite{Coin}.

If $C^{gl}$ (resp. $C^\pm$) is the unnormalized chain
complex associated to $\mathcal{N}^{gl}$ (resp. $\mathcal{N}^{\pm}$)
and if $\CR^{gl}$ (resp. $\CR^\pm$) is the chain complex which is the
quotient of $C^{gl}$ (resp. $C^\pm$) by the subcomplex generated
loosely speaking by the elements without volume (the so-called
thin elements, cf. Section~\ref{genera2}), then

\begin{conj}\label{c1}\cite{Coin} (The thin elements conjecture)
Let $K$ be a precubical set. Then the chain complex morphisms
$C^{gl}(F(K))\rightarrow \CR^{gl}(F(K))$ and $C^{\pm}(F(K))\rightarrow
\CR^{\pm}(F(K))$ induce isomorphisms in homology. \end{conj}

As a matter of fact, it is also natural to think that

\begin{conj}\cite{Coin} (The thin elements conjecture) Let $G$
be a globular set. Let $\F G$ be the free strict globular
$\omega$-category generated by $G$ (see Section~\ref{left-glob}).
Then the chain complex morphisms $C^{gl}(\F G)\rightarrow
\CR^{gl}(\F G)$ and $C^{\pm}(\F G)\rightarrow
\CR^{\pm}(\F G)$ induce isomorphisms in homology. \end{conj}

As explained in \cite{Coin}, the thin elements conjecture is closely related
to the presence or not of relations like ``$a*b=c*d$'' with $(a,b)\neq (c,d)$
in the composition laws of $\C$. Therefore it
is plausible to think that

\begin{conj}\label{c2} (The extended  thin elements conjecture)
Let $K$ be a precubical set.  Then the chain complex morphisms
\[C^{gl}(\widehat{F(K)})\rightarrow \CR^{gl}(\widehat{F(K)})\] and
\[C^{\pm}(\widehat{F(K)})\rightarrow \CR^{\pm}(\widehat{F(K)})\] induce
isomorphisms in homology. \end{conj}

Conjecture~\ref{c2} can also be formulated for free
$\omega$-categories generated by globular sets. In fact, one can
even put forward a slightly more general statement\thinspace:

\begin{conj} Le $\C$ be a strict globular non-contracting
$\omega$-category such that the chain complex morphisms
$C^{gl}(\C)\rightarrow \CR^{gl}(\C)$ and $C^{\pm}(\C)\rightarrow
\CR^{\pm}(\C)$ induce isomorphisms in homology. Then the chain complex
morphisms $C^{gl}(\widehat{\C})\rightarrow \CR^{gl}(\widehat{\C})$ and
$C^{\pm}(\widehat{\C})\rightarrow \CR^{\pm}(\widehat{\C})$ induce
isomorphisms in homology as well. \end{conj}

Chain complexes $\CF^\pm$ called the
\textit{formal branching} and \textit{the formal merging complexes}
were introduced in \cite{Coin} (cf. Definition~\ref{formalcoin}). The
\textit{formal globular complex} $\CF^{gl}$ was introduced in
\cite{sglob}\thinspace: the only difference with Definition~\ref{formalcoin} is
that the relation $x=x*_0y$ is removed.

\bd\label{new_glob}\cite{sglob}
Let $\C$ be a non-contracting $\omega$-category. If $X$ is a set, let $\Z X$
be the free abelian group generated by the set $X$. Let $tr^{n-1}\C$ be the
set of morphisms of $\C$ of dimension at most $n-1$. Let
\begin{itemize}
\item $\CF^{gl}_0(\C):=\Z\C_0\otimes\Z\C_0\iso \Z (\C_0\p\C_0)$
\item $\CF^{gl}_1(\C):=\Z\C_1$
\item $\CF^{gl}_n(\C)= \Z\C_n/\{x*_1y=x+y, \dots , x
*_{n-1}y=x+y\hbox{ mod }\Z tr^{n-1}\C\}$ for $n\geq 2$
\end{itemize}
with the differential map $s_{n-1}-t_{n-1}$ from $\CF^{gl}_n(\C)$ to
$\CF^{gl}_{n-1}(\C)$ for $n\geq 2$ and $s_0\otimes t_0$ from
$\CF^{gl}_1(\C)$ to $\CF^{gl}_0(\C)$ where $s_i$ (resp. $t_i$) is
the $i$-dimensional source (resp. target) map. This chain complex is called
the formal globular complex. The associated homology is denoted by
$\HF^{gl}(\C)$ and is called the formal globular homology of $\C$.
\ed

Using results of \cite{Coin,sglob}, we already know that the folding
operators (cf. Section~\ref{genera3}) of respectively the globular
nerve and the branching/merging semi-cubical nerves induce morphisms
of chain complexes $\CF^{gl}(\C)\rightarrow \CR^{gl}(\C)$ and
$\CF^{\pm}(\C)\rightarrow \CR^{\pm}(\C)$ which are onto for any strict
non-contracting globular $\omega$-category $\C$. The following
conjectures were then put forward\thinspace:

\begin{conj}\label{c3}\cite{Coin,sglob} The folding operators (cf. Section~\ref{genera3})
of respectively the globular nerve and the branching/merging
semi-cubical nerves induce quasi-isomorphisms of chain complexes
\[\CF^{gl}(\C)\rightarrow \CR^{gl}(\C)\] and
\[\CF^{\pm}(\C)\rightarrow \CR^{\pm}(\C)\]
for any strict non-contracting globular $\omega$-category $\C$. \end{conj}

We can now answer the question above.

\bth\label{reponse} Let $K$ be a precubical set. Assume that Conjecture~\ref{c1}, Conjecture~\ref{c2} and Conjecture~\ref{c3} above
hold. Then the morphisms of chain complexes \[C^{gl}(F(K))\rightarrow
C^{gl}(\widehat{F(K)})\] and \[C^{\pm}(F(K))\rightarrow
C^{\pm}(\widehat{F(K)})\] induce isomorphisms in homology. \eth

\bpf Let us consider the following commutative diagram\thinspace:
\[
\xymatrix{
H^{gl}(F(K))\fd{}\fr{\iso}&  \HR^{gl}(F(K)) \fd{}& \fl{\iso}\HF^{gl}(F(K))\fd{}\\
H^{gl}(\widehat{F(K)}) \fr{\iso}& \HR^{gl}(\widehat{F(K)}) & \fl{\iso}\HF^{gl}(\widehat{F(K)})}
\]
By an easy calculation, one can check that the linear map
$\HF^{gl}(F(K))\rightarrow \HF^{gl}(\widehat{F(K)})$ is an
isomorphism. Hence the result for the globular homology. The argument
is similar for both branching and merging semi-cubical homology theories.
\epf

What do we get by working with $\widehat{F(K)}$ instead of $F(K)$ ?
The globular nerve becomes Kan. Indeed its simplicial part is nothing
else but the simplicial nerve of the path $\omega$-category $\P
\widehat{F(K)}$ of $\widehat{F(K)}$ which turns out to be a strict
globular $\omega$-groupoid \cite{sglob}.

So now we can ask the question\thinspace: do the branching and merging
semi-cubical nerves $\mathcal{N}^-(\widehat{F(K)})$ and
$\mathcal{N}^+(\widehat{F(K)})$ satisfy the Kan condition as well ?
The answer is\thinspace: almost never. A counterexample is given at the end of
\cite{sglob}. Let us recall it. Consider the
$2$-source of $R(000)$ in Figure~\ref{partial_3cube} where  $R(0+0)$ is
removed. Consider both inclusion $\omega$-functors from $I^2$ to respectively
$R(-00)$ and $R(00-)$. Then the Kan condition fails because one cannot
make the sum of $R(-00)$ and $R(00-)$ since $R(0+0)$ is removed.

The purpose of this paper is to find two new augmented simplicial
nerves of non-contracting $\omega$-category $\mathcal{N}^{gl^-}(\C)$
and $\mathcal{N}^{gl^+}(\C)$ called the \textit{branching and merging
semi-globul\-ar nerves} such that the following statements hold\thinspace:
\begin{enumerate}
\item there exist natural morphisms of chain complexes
$C^{-}(\C)\rightarrow C^{gl^-}(\C)$, $\CR^{-}(\C)\rightarrow \CR^{gl^-}(\C)$,
$C^{+}(\C)\rightarrow C^{gl^+}(\C)$ and
$\CR^{+}(\C)\rightarrow \CR^{gl^+}(\C)$ for any non-contracting
$\omega$-category $\C$ such that the following squares of chain complexes
are commutative\thinspace:
\[
\xymatrix{
C^{-}(\C)\fd{}\fr{} & C^{gl^-}(\C)\fd{} &C^{+}(\C)\fd{}\fr{} & C^{gl^+}(\C)\fd{}\\
\CR^{-}(\C)\fr{} & \CR^{gl^-}(\C)&\CR^{+}(\C)\fr{} & \CR^{gl^+}(\C)
}\]
\item If $\C$ is an object of $\omega Cat_1^{Kan}$, then the simplicial sets
$\mathcal{N}^{gl^-}(\C)$ and $\mathcal{N}^{gl^+}(\C)$ satisfy the Kan
condition.
\item The new simplicial nerves $\mathcal{N}^{gl^-}$ and $\mathcal{N}^{gl^+}$
fit the formalism of regular cut, as the old branching and merging
nerves do. As a consequence, there exist folding operators
$\Phi^{gl^\pm}$ and $\square_n^{gl^\pm}$ inducing morphisms of chain
complexes $\CF^{gl^\pm}(\C)=\CF^\pm(\C)\rightarrow \CF^{gl^\pm}(\C)$
which are onto (and conjecturally quasi-isomorphisms) for any strict
non-contracting $\omega$-category $\C$.
\item There exist two morphisms of augmented simplicial sets
$\mathcal{N}^{gl}\rightarrow \mathcal{N}^{gl^-}$ and
$\mathcal{N}^{gl}\rightarrow \mathcal{N}^{gl^+}$ which do what we
want, that is associating in homology to any empty oriented globe its
corresponding branching and merging areas of execution paths (as in
\cite{sglob})\thinspace: cf. Figure~\ref{recap} in Section~\ref{comparison}.
\end{enumerate}

Let $K$ be a precubical set. Then there will exist a commutative diagram
\[
\xymatrix{
H^{gl^\pm}(F(K))\fd{}\fr{\iso}&  \HR^{gl^\pm}(F(K)) \fd{}& \fl{\iso}\HF^{\pm}(F(K))\fd{}\\
H^{gl^\pm}(\widehat{F(K)}) \fr{\iso}& \HR^{gl^\pm}(\widehat{F(K)}) & \fl{\iso}\HF^{\pm}(\widehat{F(K)})
}
\]
where the isomorphisms are conjectural. Notice that the formal branching (resp.
merging) homology associated to the branching (resp. merging) semi-globular
nerve is the same as that associated to the branching (resp. merging)
semi-cubical
nerve. In particular, for a strict globular $\omega$-category freely generated
by a precubical set, the branching (resp. merging) semi-globular and semi-cubical
homologies conjecturally coincide.

The reader may find that a lot of conjectures remain to be
proved. This is certainly true but they require new ideas to be
resolved and the computer-scientific motivations are helpless. The
study of HDA provides new combinatorial conjectures (exactly as in
\cite{Coin,sglob}) which need new ideas to be proved.

How to construct the branching semi-globular nerve
$\mathcal{N}^{gl^-}(\C)$ (and by symmetry the merging semi-globular
nerve $\mathcal{N}^{gl^+}(\C)$) ? We already said that the formal
homology theory associated to the branching semi-globular nerve is
expected to be the same as that of the branching semi-cubical nerve.
And the formal branching complex of a given $\omega$-category $\C$ is
exactly the formal globular complex of an $\omega$-category
``$\C/(x=x*_0y)$'' where ``$\C/(x=x*_0y)$'' would be an
$\omega$-category associated to $\C$ in which the relation $x=x*_0y$
is forced (for $dim(x)\geq 1$).

So loosely speaking, the branching semi-globular nerve of a strict
globular $\omega$-category $\C$ will be the globular nerve of
``$\C/(x=x*_0y)$''. Looking back to our counterexample, one sees that
$R(-00)$ and $R(00-)$ become composable because in ``$\C/(x=x*_0y)$'',
$t_1 R(-00) = s_1 R(00-) = R(-0-)$. The construction of
``$\C/(x=x*_0y)$'' is concretely implemented in
Section~\ref{left-glob}. The corresponding globular nerve will be
called \textit{the branching semi-globular cut}
$\mathcal{N}^{gl^-}(\C)$ of $\C$.

The organization of the paper is as follows.
Section~\ref{convention} recalls our conventions of notations and
also Steiner's formulas in $\omega$-complexes\thinspace: these
formulas are indeed crucial for several proofs of this paper. In
Section~\ref{left-glob}, the construction of ``$\C/(x=x*_0y)$''
(it is called the \textit{negative semi-path $\omega$-category}
$\P^-\C$ of $\C$) is presented.  Some elementary facts about the
latter are proved. In Section~\ref{case-hypercube}, the technical
core of the paper, the negative semi-path $\omega$-category of
the hypercube is completely calculated in any dimension (see
Theorem~\ref{ortho_cube} and Corollary~\ref{ortho_complet}). This
calculation will indeed be fundamental to construct in
Section~\ref{comparison} the canonical natural transformation
from the branching semi-cubical cut to the branching
semi-globular cut. Section~\ref{genera1} exposes a generalization
of the notion of \textit{simplicial cut} initially introduced in
\cite{sglob}.  The reason of this generalization is that the
branching semi-globular cut does not exactly match the old
definition. Section~\ref{genera2} then recalls the notion of
\textit{thin elements} and at last Section~\ref{genera3} the
notion of \textit{regular cut}. The latter definition is
necessary to construct the \textit{folding operators}.  Each
definition introduced in Section~\ref{genera1},
Section~\ref{genera2} and Section~\ref{genera3} is illustrated
with the case of the branching semi-cubical situation (i.e. the
old definition of the branching complex).  At last the definition
of the branching semi-globular nerve in Section~\ref{definition}
and the comparison with the branching semi-cubical nerve in
Section~\ref{comparison}.

Of course all constructions of this paper may be applied to the case
of the merging nerve in an obvious way. The Kan version of the merging
nerve will be called the \textit{merging semi-globular nerve}.  Figure~\ref{recap} is
a recapitulation of all simplicial constructions obtained so far
(including these of this paper).

This work is part of a research project which aims at setting up an
appropriate algebraic setting for the study deformations of HDA
which leave invariant their computer-scientific properties. See
\cite{ConcuToAlgTopo} for a sketch of a description of the project.

Let us conclude by some remarks about the terminology. A lot of
functors have been introduced so far and some coherence in their
naming is necessary. Let $\C$ be a non-contracting globular $\omega$-category.
\begin{enumerate}
\item $\P\C$ is the \textit{path $\omega$-category} of $\C$.
\item $\P^-\C$ is the \textit{negative semi-path $\omega$-category} of $\C$ and
$\P^+\C$ its \textit{positive semi-path $\omega$-category}.
\item For $\alpha\in\{-,+\}$, $\P^\alpha\C$ is the \textit{semi-path $\omega$-category}
of $\C$.
\item The nerve $\mathcal{N}^-(\C)$ is the \textit{branching semi-cubical nerve} of
$\C$ and $\mathcal{N}^+(\C)$ its \textit{merging semi-cubical nerve}.
Without further precision about $\alpha\in\{-,+\}$,  $\mathcal{N}^\alpha(\C)$ is the
\textit{semi-cubical nerve}. It is called \textit{corner nerve} in previous publications.
\item The simplicial homology shifted by one of the branching semi-cubical nerve is
called the \textit{branching semi-cubical homology} and the simplicial homology shifted by one of the merging semi-cubical nerve is
called the \textit{merging semi-cubical homology}.
\item The nerve $\mathcal{N}^{gl^-}(\C)$ is the \textit{branching semi-globular nerve} of
$\C$ and $\mathcal{N}^{gl^+}(\C)$ its \textit{merging semi-globular nerve}. Without further precision about $\alpha\in\{-,+\}$,  $\mathcal{N}^{gl^\alpha}(\C)$ is the
\textit{semi-globular nerve}.
\item The simplicial homology shifted by one of the branching semi-globular nerve is
called the \textit{branching semi-globular homology} and is
denoted by $H_{n+1}^{gl^-}(\C):=H_n(\mathcal{N}^{gl^-}(\C))$\thinspace;
the simplicial homology shifted by one of the merging
semi-globular nerve is called the \textit{merging semi-globular
homology} and is denoted by
$H_{n+1}^{gl^+}(\C):=H_n(\mathcal{N}^{gl^+}(\C))$.
\end{enumerate}

\section{Preliminaries}\label{convention}

The reader who is familiar with papers \cite{Gau,Coin,sglob} may
want to skip this section.

\bd\label{omega_categories} \cite{Brown-Higgins0,oriental,Tensor_product}
An \textit{$\omega$-category} is a set $A$
endowed with two families of maps $(s_n=d_n^-)_{n\geqslant 0}$ and
$(t_n=d_n^+)_{n\geqslant 0}$ from $A$ to $A$ and with a family of partially
defined 2-ary operations $(*_n)_{n\geqslant 0}$ where for any
$n\geqslant 0$, $*_n$ is a map from $\{(a,b)\in A\p A,
t_n(a)=s_n(b)\}$ to $A$ ($(a,b)$ being carried over $a *_n b$) which
satisfies the following axioms for all $\alpha$ and $\beta$ in
$\{-,+\}$\thinspace:

\begin{enumerate}
\item $d_m^\beta d_n^\alpha x=
\left\{\begin{CD}d_m^\beta x \hbox{  if $m<n$}\\  d_n^\alpha x \hbox{  if $m\geqslant n$}
\end{CD}\right.$
\item $s_n x *_n x= x *_n t_n x = x$
\item if  $x *_n y$ is well-defined, then  $s_n(x *_n y)=s_n x$, $t_n(x *_n y)=t_n y$
and for  $m\neq n$, $d_m^\alpha(x *_n y)=d_m^\alpha x *_n d_m^\alpha y$
\item as soon as the two members of the following equality exist, then
$(x *_n y) *_n z= x *_n (y *_n z)$
\item if $m\neq n$ and if the two members of the equality make sense, then
$(x *_n y)*_m (z*_n w)=(x *_m z) *_n (y *_m w)$
\item for any  $x$ in $A$, there exists a natural number $n$ such that $s_n x=t_n x=x$
(the smallest of these numbers is called the dimension of $x$ and is denoted
by $dim(x)$).
\end{enumerate}
\ed

We will sometimes use the notations $d^-_n:=s_n$ and $d^+_n=t_n$. If
$x$ is a morphism of an $\omega$-category $\C$, we call $s_n(x)$ the
$n$-source of $x$ and $t_n(x)$ the $n$-target of $x$.  The category of
all $\omega$-categories (with the obvious morphisms) is denoted by
$\omega Cat$. The corresponding morphisms are called
$\omega$-functors.

Roughly speaking, our $\omega$-categories are strict, globular and
contain only morphisms of finite dimension.  If $\C$ is an object of
$\omega Cat$, then $tr^n\C$ denotes the set of morphisms of $\C$ of
dimension at most $n$ and $\C_n$ the set of morphisms of dimension
exactly $n$.  The map $tr^n$ induces for any $n\geq 0$ a functor from
$\omega Cat$ to itself.

\bd Let $\C$ be an $\omega$-category. Then $\C$ is said
\textit{non-contracting} if and only if for any morphism $x$ of
$\C$ of dimension at least $1$, then $s_1x$ and $t_1x$ are
$1$-dimensional. \ed

\bd\label{noncontractant} Let $f$ be an $\omega$-functor from $\C$ to
$\D$.  The morphism $f$ is \textit{non-contract\-ing} if for any
$1$-dimensional $x\in \C$, the morphism $f(x)$ is a $1$-dimensional
morphism of $\D$ (a priori, $f(x)$ could be either $0$-dimensional
or $1$-dimensional).  \ed

\begin{nota}\label{catnoncontract} The category of non-contracting
$\omega$-categories with the non-contracting $\omega$-functors is
denoted by $\omega Cat_1$.  \end{nota}

\begin{nota} If $\C$ is a non-contracting
$\omega$-category, one can consider the $\omega$-category $\P\C$
such that $(\P\C)_i:=\C_{i+1}$, $*_i^{\P\C}=*_{i+1}^{\C}$,
$*_i^{\P\C}=*_{i+1}^{\C}$ and $*_i^{\P\C}=*_{i+1}^{\C}$ for any
$i\geq 0$.  The $n$-source (resp. the $n$-target, the
$n$-dimensional composition law) of $\P\C$ will still be denoted
by $s_{n+1}$ (resp. $t_{n+1}$, $*_{n+1}$) to avoid possible
confusion. \end{nota}

\bd\label{path-cat}  The $\omega$-category $\P\C$ is called the
\textit{path $\omega$-category} of $\C$. The mapping $\P$ yields a
functor from $\omega Cat_1$ to $\omega Cat$. \ed

Fundamental examples of $\omega$-categories are the
$\omega$-category associated to the $n$-dimens\-ional cube and
the one associated to the $n$-dimensional simplex (the former is
denoted by $I^n$ and the latter by $\Delta^n$).  Both families of
$\omega$-categories can be characterized in the same way. The
first step consists of labeling all faces of the $n$-cube and of
the $n$-simplex. For the $n$-cube, this consists of considering
all words of length $n$ in the alphabet $\{-,0,+\}$, one word
corresponding to the barycenter of a face (with $00\dots 0 \hbox{
($n$ times)}=:0_n$ corresponding to its interior). As for the
$n$-simplex, its faces are in bijection with strictly increasing
sequences of elements of $\{0,1,\dots,n\}$. A sequence of length
$p+1$ will be of dimension $p$. If $x$ is a face, let $R(x)$ be
the set of subfaces of $x$ (including $x$) seen respectively as a
sub-cube or a sub-simplex. If $X$ is a set of faces, then let
$R(X)=\bigcup_{x\in X}R(x)$. Notice that $R(X\cup Y)=R(X)\cup
R(Y)$, $R(X\cap Y)=R(X)\cap R(Y)$ and $R(\{x\})=R(x)$. Then $I^n$
and $\Delta^n$ are the unique $\omega$-categories such that the
underlying set is a subset of $\{R(X),X \hbox{ set of faces}\}$
satisfying the following properties\thinspace:

\begin{enumerate}
\item  For $x$ a $p$-dimensional face of $I^n$ (resp. $\Delta^n$) with $p\geq 1$,
$s_{p-1}(R(x))=R(s_x)$ and $t_{p-1}(R(x))=R(t_x)$ where $s_x$ and
$t_x$ are the sets of faces defined below.
\item If $X$ and $Y$ are two elements of $I^n$ (resp. $\Delta^n$) such that
$t_p(X)=s_p(Y)$ for some $p$, then $X\cup Y$ belongs to $I^n$
(resp. $\Delta^n$) and $X\cup Y=X *_p Y$.
\end{enumerate}

Only the definitions of $s_x$ and $t_x$ differ between the
construction of the family of $\Delta^n$ and that of $I^n$. Let
us give the computation rule in some examples. For the cube, the
$i$-th zero is replaced by $(-)^i$ (resp. $(-)^{i+1}$) for $s_x$
(resp.  $t_x$). For example, one has
$s_{0+00}=\{\hbox{-+00},\hbox{0++0},\hbox{0+0-}\}$ and
$t_{0+00}=\{\hbox{++00},\hbox{0+-0},\hbox{0+0+}\}$. For the
simplex, $s_{(04589)}=\{(4589),(0489),(0458)\}$ (the elements in
odd position are removed) and $t_{(04589)}=\{(0589),(0459)\}$
(the elements in even position are removed).

One can also notice that if $X$ is an element of $I^n$ or
$\Delta^n$, then $X=R(X)$.

Both $\Delta^n$ and $I^n$ are examples of $\omega$-complexes in
the sense of \cite{MR99e:18008} where the atoms are the elements
of the form $R(\{x\})$ where $x$ is a face of $\Delta^n$ (resp.
$I^n$). In such situations, Steiner's paper \cite{MR99e:18008}
proves that the calculation rules are very simple and that they
can be summarized as follows\thinspace:

\begin{enumerate}
\item If $X$ and $Y$ are two elements of $I^n$ (resp. $\Delta^n$) such that
$t_p(X)=s_p(Y)$ for some $p$, then not only $X\cup Y$ belongs to $I^n$
(resp. $\Delta^n$) and $X\cup Y=X *_p Y$ but also $t_p(X)=s_p(Y)=X\cap Y$.
\item For $x$ a $p$-dimensional face of $I^n$ (resp. $\Delta^n$) with $p\geq 0$,
set $\de^-R(x):= s_{p-1}R(x)=R(s_x)$ if $p>0$ and $\de^-R(x):=\emptyset$ if $p=0$ and
set $\de^+R(x):= t_{p-1}R(x)=R(t_x)$ if $p>0$ and $\de^+R(x):=\emptyset$ if $p=0$\thinspace;
then for any $X$ in $I^n$ (resp. $\Delta^n$), then
$$d_n^\alpha X = \left(\bigcup_{a\in X, dim(a)\leq n} R(a)\right)\backslash \left(\bigcup_{b\in X, dim(b)=n+1}\left(R(b)\backslash \de^\alpha R(b)\right)\right) $$
with $d_n^-=s_n$ and $d_n^+=t_n$.
\end{enumerate}

In a simplicial set, the face maps are always denoted by $\de_i$, the
degeneracy maps by $\epsilon_i$. Here are the other conventions about
simplicial sets and simplicial homology theories
(see for example \cite{May} for further information)\thinspace:

{
\begin{enumerate}
\item $Sets$\thinspace: category of sets
\item $Sets^{\Delta^{op}}$\thinspace: category of simplicial sets
\item $Sets^{\Delta^{op}}_+$\thinspace: category of augmented simplicial sets
\item $\comp$\thinspace: category of chain complexes of abelian groups
\item $C(A)$\thinspace: unnormalized chain complex of the simplicial set $A$
\item $H_*(A)$\thinspace: simplicial homology of a simplicial set $A$
\item $Ab$\thinspace: category of abelian groups
\item $\id$ \thinspace: identity map
\item $\Z S$\thinspace: free abelian group generated by the set $S$
\end{enumerate}}

\section{The negative semi-path $\omega$-category of a non-contracting $\omega$-category}\label{left-glob}

Before going further in the construction of the negative semi-path
$\omega$-category of a non-contracting $\omega$-category, one
needs to recall some well-known facts about globular sets
\cite{petittoposglob}.  Let us consider the small category $Glob$
defined as follows\thinspace: the objects are all natural numbers
and the arrows are generated by $s$ and $t$ in $Glob(m,m-1)$ for
any $m>0$ and by the relations $s\circ s=s\circ t$, $t\circ
s=t\circ t$.  By definition, a globular set is a covariant
functor from $Glob$ to the category of sets $Sets$. Let us denote
by $\mathcal{G}_\omega$ the corresponding category. Let $\U$ be
the forgetful functor from $\omega Cat$ to $\mathcal{G}_\omega$.
One can prove by standard categorical arguments the existence of
a left adjoint $\F$ for $\U$ (see \cite{MR2000c:18006} for an
explicit construction of this left adjoint).

\bd If $G$ is a globular set, then the $\omega$-category $\F G$
is called the free $\omega$-category generated by the globular set
$G$. \ed

Let $\C$ be a non-contracting $\omega$-category. Let us denote by $\Rm$ the
equivalence relation which is the reflexive, symmetric and transitive closure of
the subset \[\left\{(x,x *_0 y); (x,y, x*_0y)\in \P\C\p\P\C\p\P\C\right\}\]
of $\P\C\p \P\C$.
Now consider the underlying globular set $\U\P\C$ of the path
$\omega$-category of a non-contracting $\omega$-category $\C$.  Since
$\C$ is non-contracting, $s_i(\P\C)\subset \P\C$ and $t_i(\P\C)\subset
\P\C$ for any $i\geq 1$ and the maps $s_i$ and $t_i$ from $\P\C$ to
itself pass to the quotient for $i\geq 1$ because $s_i(x*_0 y)=s_i x
*_0 s_i y$ and $t_i(x*_0 y)=t_i x *_0 t_i y$ for any $i\geq 1$.

Let $\phi:\U\P\C\rightarrow \U\P\C/\Rm$ be the canonical morphism of
globular sets induced by $\Rm$.  The identity map $\U \P\C\rightarrow
\U \P\C$ provides the canonical morphism of $\omega$-categories
$\eta_{\P\C}:\F \U \P\C \rightarrow \P\C$, $\eta$ being the counit of
the adjunction $(\F,\U)$. Then let us consider the following pushout
in $\omega Cat$\thinspace:
\[
\xymatrix{
{\F \left(\U \P\C\right)} \ar@{->}[r]^-{\F(\phi)} \fd{\eta_{\P\C}} & {\F \left(\U \P\C/\Rm\right)} \ar@{->}[d]^{A}   \\
{\P\C} \fr{h^-}  &   \cocartesien {\P^-\C}
}
\]

\bd Let $\C$ be a non-contracting $\omega$-category. The $\omega$-category
$\P^-\C$ defined above is called the negative semi-path $\omega$-category
of $\C$. \ed

The negative semi-path $\omega$-category $\P^-\C$ of $\C$ intuitively
contains the germs of non-constant execution paths of $\C$ beginning
in the same way and the germs of higher dimensional homotopies between
them. The construction above yields a functor $\P^-:\omega
Cat_1\rightarrow \omega Cat$ and a natural transformation
$h^-:\P\rightarrow \P^-$ between functors from $\omega Cat_1$ to
$\omega Cat$.

\bp\label{sur} Consider the construction of $\P^-\C$ above. Then any element of
$\P^-\C$ is a composite of elements of the form $A(x)$ with
$x\in {\F \left(\U \P\C/\Rm\right)}$. \ep

\bpf First let us make a short digression. Let $\xymatrix@1{{\C}
\ar@<1ex>@{->}[r]^{f}\ar@<-1ex>@{->}[r]_{g}& {\D}}$ be a pair of
$\omega$-functors $f$ and $g$ from an $\omega$-category $\C$ to
an $\omega$-category $\D$. Then the coequalizer $h:\D\rightarrow
\mathcal{E}$ of $f$ and $g$ always exists in $\omega Cat$ and any
element of $\mathcal{E}$ is a composite of elements of the form
$h(x)$ where $x$ runs over $\D$ (otherwise take the image of $\D$
in the coequalizer\thinspace: this image still satisfies the
universal property of the equalizer).

Since the canonical $\omega$-functor $h^-\oplus A:\P\C\oplus \F
\left(\U \P\C/\Rm\right)\rightarrow \P^-\C$ is the coequalizer of
$\F(\phi)$ and $\eta_{\P\C}$, the previous remark does apply (the
symbol $\oplus$ meaning the direct sum in $\omega Cat$, which
coincides with the disjoint union). But $h^-\oplus A$ identifies
any element of $\P\C$ with an element of $\F \left(\U
\P\C/\Rm\right)$ therefore any element of $\P^-\C$ is a composite
of elements of the form $A(x)$ with $x$ running over $\F \left(\U
\P\C/\Rm\right)$. \epf

\bth (Universal property satisfied by $\P^-$) Let $\C$ be a non-contracting
$\omega$-categ\-ory. Let $\D$ be an object of $\omega Cat$. Let
$\mu:\P\C\rightarrow \D$ be an $\omega$-functor such that for any
$x,y\in \P\C$, $x\Rm y$ implies $\mu(x)=\mu(y)$ in $\D$. Then there
exists a unique $\omega$-functor $\overline{\mu}:\P^-\C\rightarrow \D$
such that $\mu=\overline{\mu}\circ h^-$.
\eth

\bpf Let $\mu:\P\C\rightarrow \D$ be an $\omega$-functor such that for
any $x,y\in \P\C$, $x\Rm y$ implies $\mu(x)=\mu(y)$ in $\D$. Then
$\mu$ induces a morphism of globular sets
$\U(\mu):\U\P\C\rightarrow \U \D$ and by hypothesis, $\U(\mu)$ gives
rise to a morphism of globular
sets $\overline{\U(\mu)}:\U\P\C/\Rm \rightarrow \U \D$ such that
$\overline{\U(\mu)}\circ \phi = \U(\mu)$. Then the composite
\[
\xymatrix{\F(\U\P\C/\Rm)\ar@{->}[rr]^{\F(\overline{\U(\mu)})}&& \F\U \D \fr{\eta_\D}&\D}
\]
yields an $\omega$-functor from $\F(\U\P\C/\Rm)$ to $\D$. Since
$\eta:\F\U\rightarrow Id_{\omega Cat}$ is a natural
transformation, one gets the commutative diagram
\[
\xymatrix{
\F\U(\P\C)\fr{\F\U(\mu)} \fd{\eta_{\P\C}}& \F\U(\D)\fd{\eta_{\D}}\\
\P\C \fr{\mu}  & \D }
\]
so the equality $\eta_\D \circ \F\U(\mu)=\mu \circ \eta_{\P\C}$ holds.
Therefore $\eta_\D \circ \F(\overline{\U(\mu)}) \circ \F(\phi)= \eta_\D \circ \F\U(\mu)= \mu \circ \eta_{\P\C}$. One then obtains the commutative diagram
\[
\xymatrix{
{\F \left(\U \P\C\right)} \ar@{->}[r]^-{\F(\phi)} \fd{\eta_{\P\C}} & {\F \left(\U \P\C/\Rm\right)} \ar@{->}[d]^{A}   \ar@/^15pt/[rdd]^{\eta_\D \circ \F(\overline{\U(\mu)})}&\\
{\P\C} \fr{h^-}  \ar@/_15pt/[rrd]_{\mu}&   \cocartesien {\P^-\C}  &\\
&& \mathcal{D}
}
\]
Therefore there exists a unique natural transformation
$\overline{\mu}:\P^-\C\rightarrow \D$ such that
$\overline{\mu}\circ A = \eta_\D\circ \F(\overline{\U(\mu)})$ and
such that $\overline{\mu}\circ h^-= \mu$, i.e. making the following diagram commutative\thinspace:
\[
\xymatrix{
{\F \left(\U \P\C\right)} \ar@{->}[r]^-{\F(\phi)} \fd{\eta_{\P\C}} & {\F \left(\U \P\C/\Rm\right)} \ar@{->}[d]^{A} \ar@{->}[rd]^{\F(\overline{\U(\mu)})} &  \\
{\P\C} \fr{h^-} \ar@{->}@/_20pt/[rrd]_{\mu} &   \cocartesien {\P^-\C} \ar@{->}[rd]^{\overline{\mu}} & {\F\U \D} \fd{\eta_\D}\\
&& {\D}
}
\]

Suppose that there exists another $\nu:\P^-\C\rightarrow \D$ such that
$\nu \circ h^-= \mu$. Proving that $\nu=\overline{\mu}$ is equivalent
to proving that $\nu \circ A = \overline{\mu}\circ A$ by
Proposition~\ref{sur}.  So one is reduced to checking the equality
$\nu \circ A =\eta_\D \circ \F(\overline{\U(\mu)})$.  The $\omega$-functor
$\F(\phi)$ is clearly surjective on the underlying sets. Therefore
proving $\nu \circ A = \eta_\D \circ \F(\overline{\U(\mu)})$ is
equivalent to proving $\nu \circ A \circ \F(\phi) = \eta_\D \circ
\F(\overline{\U(\mu)})\circ \F(\phi)$. But one has $\nu \circ A \circ
\F(\phi)=\nu \circ h^- \circ \eta_{\P\C}=\mu \circ \eta_{\P\C}=
\eta_\D \circ \F(\overline{\U(\mu)}) \circ \F(\phi)$, which concludes the proof.
\epf

By convention, in any of the above $\omega$-categories arising from
$\C$, the $n$-source (resp. the $n$-target, the $n$-dimensional
composition law) will be still denoted by $s_{n+1}$, (resp. $t_{n+1}$,
$*_{n+1}$), like for $\P\C$.  The calculation rules in $\P^-\C$ are
summarized in the following theorem\thinspace:

\bth\label{calc_rules} (Calculation rules in $\P^-\C$)
Let $\C$ be a non-contracting $\omega$-category.
Then any element of $\P^-\C$ is a composite of elements of the form
$h^-(x)$. Moreover if $x$ and $y$ are two elements of $\P\C$ such that
$x *_p y$ exists in $\P\C$ for some $p\geq 1$, then $h^-(x *_p
y)=h^-(x) *_p h^-(y)$. And if $x$ and $y$ are two elements of $\P\C$
such that $x *_0 y$ exists in $\P\C$, then $h^-(x*_0y)=h^-(x)$.
\eth

\bpf By Proposition~\ref{sur}, any element of $\P^-\C$ is a
composite of elements of the form $A(x)$ with $x\in {\F \left(\U
\P\C/\Rm\right)}$.  Since $\F(\phi)$ is clearly surjective on the
underlying sets, any element of $\P^-\C$ is a composite of
elements of the form $A\circ \F(\phi)(x)=h^-(\eta_{\P\C}(x))$ with
$x\in {\F \left(\U\P\C\right)}$. So any element of $\P^-\C$ is a
composite of elements of the form $h^-(x)$ with $x$ running over
$\P \C$.  The last part of the statement of the theorem is a
consequence of the fact that $h^-$ is an $\omega$-functor and of
the universal property satisfied by $\P^-\C$. \epf

Loosely speaking, the $\omega$-category $\P^-\C$ is the quotient
of the free $\omega$-category generated by the equivalence classes
of $\Rm$ in $\U\P\C$ by the calculation rules of $\C$. The
calculation rules in $\P^-\C$ are more explicitly described as
follows.
\begin{enumerate}
\item If $x$ and $y$ are two morphisms of $\P\C$ such that
$x\Rm y$, then $\phi(x)=\phi(y)$ and therefore
$x$ and $y$ give rise to the same element in $\P^-\C$\thinspace: in other
terms $h^-(x)=h^-(y)$.
\item If $x$ and $y$ are two morphisms of $\P\C$ such that
$t_p x\Rm s_p y$ for some $p\geq 1$, then
$\phi(t_p x)=\phi(s_p y)$ and
the corresponding elements can be composed in $\P^-\C$. Therefore
$h^-(x) *_p h^-(y)$ exists although $x*_p y$ does not necessarily
exist in $\P\C$.
\item If moreover $x$ and $y$ are two elements of $\P\C$ such that
this time $t_p x= s_p y$ for some $p\geq 1$, then
$h^-(x) *_p h^-(y)$ is the image of $x *_p y$ by $h^-$, and therefore
$h^-(x *_p y)=h^-(x) *_p h^-(y)$.
\item If $x$ and $y$ are two morphisms of $\P\C$ such that $x*_0y\in \P\C$,
then $h^-(x*_0y)=h^-(x)$.
\end{enumerate}

\bp\label{grading}
Let $\C$ be a non-contracting $\omega$-category. The $0$-source map
$s_0$ of $\C$ induces $\C_0$-gradings on $\F\U\P\C$,
$\F(\U\P\C/\mathcal{R}^-)$, $\P\C$ and therefore on $\P^-\C$ as well.
Let us denote by $G^{\alpha,-}\F\U\P\C$,
$G^{\alpha,-}\F(\U\P\C/\mathcal{R}^-)$, $G^{\alpha,-}\P\C$ and
$\P^-_\alpha\C$ the fiber of $s_0$ in the respective $\omega$-categories
over $\alpha\in \C_0$. Then one has the pushout
\[
\xymatrix{
{G^{\alpha,-}\F \left(\U \P\C\right)} \fr{} \fd{} & {G^{\alpha,-}\F \left(\U \P\C/\Rm\right)}\ar@{->}[d]\\
{G^{\alpha,-}\P\C} \fr{} & {\P^-_\alpha\C} \cocartesien
}
\]
\ep

\bpf This is due to the fact that $s_0(x *_0 y) = s_0 x$ for any
$x,y\in \P\C$.  \epf

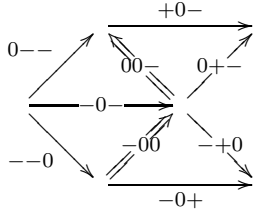
\begin{figure}
\[
\xymatrix{
&\ar@{->}[rr]^{+0-}& &\\
\ar@{->}[ru]^{0--}\ar@{->}[rr]|{-0-}\ar@{->}[rd]_{--0}&&\ar@{->}[ru]|{0+-}\ar@{->}[dr]|{-+0}\ff{lu}{00-}&&\\
&\ar@{->}[rr]_{-0+}\ff{ru}{-00}&&\\
}
\]
\caption{Part of the $2$-source of the $3$-cube}
\label{partial_3cube}
\end{figure}

As an example, consider the $\omega$-category $\C$ of
Figure~\ref{partial_3cube} which is a part of the $2$-source side
of the $3$-cube. The underlying set of the $\omega$-category
$\P^-_{R(---)}\C$ is equal to {\scriptsize $$ \left\{h^-(R(--0)),
h^-(R(-0-)), h^-(R(0--)), h^-(R(-00)),h^-(R(00-)),
h^-(R(-00))*_1h^-(R(00-))\right\}
$$}
Notice that $h^-(R(-00))$ and $h^-(R(00-))$ become composable in $\P^-\C$, although
they are not composable in the initial $\omega$-category $\C$.

Now let us recall the notion of strict globular $\omega$-groupoid\thinspace:

\bd \cite{Brown_cube} \label{def_groupoid}
Let $\C$ be a strict globular $\omega$-category. Then $\C$ is a
strict globular $\omega$-groupoid if and only if for any $p$-morphism $A$ of $\C$
with $p\geq 1$ and any $r\geq 0$, then there exists $A'$ (a priori depending on $r$)
such that $A*_r A'=s_r A=t_r A'$ and
$A'*_r A=s_r A'=t_rA$. \ed

\bth\label{groupoid_aussi} Let $\C$ be a  non-contracting $\omega$-category. If
$\P\C$ is a strict globular $\omega$-groupoid, then
$\P^-\C$ is a strict globular $\omega$-groupoid as well.
\eth

\bpf Any element $X$ of $\P^-\C$ is a composite of elements $h^-(x_i)$
for $i=1,\dots,n$ by Theorem~\ref{calc_rules}. For a given $X$,
let us call the smallest possible $n$ the length of $X$. Now we check
by induction on $n$ the property $P(n)$\thinspace: ``for any $X\in \P^-\C$ of
length at most $n$, for any $r\geq 1$, there exists $Y$ such that $X
*_r Y=s_rX=t_rY$ and $Y*_r X=s_rY=t_rX$''. For $n=1$, this is an
immediate consequence of the fact that $\P\C$ is an $\omega$-groupoid.
Now suppose $P(n)$ proved for $n=n_0$ with $n_0\geq 1$. Let $X$ be an
element of $\P^-\C$ of length $n_0+1$. Then $X= X_1 *_p X_2$ for some
$p\geq 1$ and with $X_1$ and $X_2$ of length at most $n_0$. Let $r\geq
1$. Let $Y_1$ (resp.  $Y_2$) be an inverse of $X_1$ (resp. $X_2$) for
$*_r$. If $r\neq p$, then $(X_1 *_p X_2) *_r (Y_1 *_p Y_2)= s_r X_1
*_p s_r X_2 = s_r (X_1 *_p X_2)$ so $Y_1 *_p Y_2$ is a solution. If
$r=p$, then $(X_1 *_p X_2) *_p (Y_2 *_p Y_1)= s_p X_1$, so $Y_2 *_p Y_1$ is now
a solution.  \epf

\bd For $\C$ a non-contracting $\omega$-category, the $\omega$-category
$\P^-\C$ is called the \textit{negative} semi-path $\omega$-category of $\C$.
\ed

In the sequel, all non-contracting $\omega$-categories $\C$ will be
supposed to have a path $\omega$-category $\P\C$ which is a strict
globular $\omega$-groupoid.

\section{The semi-path $\omega$-category of the hypercube}\label{case-hypercube}

Some $\omega$-functors will be constructed in this section using the
classical tool of filling of shells.

\bd In a simplicial set $A$, a $n$-shell is a family $(x_i)_{i=0,\dots,n+1}$
of $(n+2)$ $n$-simplexes of $A$ such that for any $0\leq i<j\leq n+1$,
$\de_i x_j=\de_{j-1} x_i$. \ed

\bp\label{filling_simp} \cite{sglob}
Let $\C$ be a non-contracting $\omega$-category. Consider
a $n$-shell \[(x_i)_{i=0,\dots,n+1}\] of the globular simplicial nerve
of $\C$. Then
\begin{enumerate}
\item The labeling defined by $(x_i)_{i=0,\dots,n+1}$ yields an
$\omega$-functor $x$ (and necessarily exactly one)
from $\Delta^{n+1}\backslash \{(01\dots n+1)\}$ to $\P\C$.
\item Let $u$ be a morphism of $\C$ such that
\[s_{n}u=x\left(s_n R((01\dots n+1))\right)\]
and
\[t_{n}u=x\left(t_n R((01\dots n+1))\right)\]
Then there exists one and only
one $\omega$-functor still denoted by $x$ from $\Delta^{n+1}$ to $\P\C$ such that
for any $0\leq i\leq n+1$, $\de_i x = x_i$ and $x((01\dots n+1))=u$.
\end{enumerate}
\ep

If $(\sigma_0<\dots <\sigma_r)$ is a face of $\Delta^n$, then
let $$\phi_n^-(((\sigma_0<\dots <\sigma_r))=k_1\dots k_{n+1}$$
where $k_{\sigma_i +1}=0$ and the other $k_j$ are equal to $-$. For
$n=2$, the set map $\phi$ looks as follows\thinspace:
\beas
&& \phi_2^-\thinspace:  (012) \mapsto 000 \\
&& \phi_2^-\thinspace: (01)  \mapsto 00- \\
&& \phi_2^-\thinspace: (02) \mapsto 0-0 \\
&& \phi_2^-\thinspace: (12) \mapsto -00 \\
&& \phi_2^-\thinspace: (0)  \mapsto 0-- \\
&& \phi_2^-\thinspace: (1)  \mapsto -0- \\
&& \phi_2^-\thinspace: (2)  \mapsto --0
\eeas

If $x$ is a face of $I^{n+1}$ which belongs to the image of
$\phi_n^-$, then denote by $x^\bot$ the unique face of $\Delta^n$ such
that $\phi_n^-(x^\bot)=x$. Notice that for any face $y$ of $\Delta^n$,
then $\phi_n^-(y)^\bot = y$.

\bp
\label{sous-delta}
Let us denote by $\Delta^n_i$ for $i=0,\dots,n$ and $n\geq 1$
the $\omega$-subcategory
of $\Delta^n$ obtained by keeping only the strictly increasing sequences
$(\sigma_0<\dots <\sigma_r)$ of $\{0,1,\dots,n\}$ such that
$i\notin \{\sigma_0,\dots,\sigma_r\}$. Then  the set map
$$(\sigma_0<\dots<\sigma_\ell<\widehat{i}\leq \sigma_{\ell+1}<\dots \sigma_r)\mapsto
(\sigma_0<\dots<\sigma_\ell<\sigma_{\ell+1}+1<\dots \sigma_r+1)
$$
(the notation
$(\sigma_0<\dots<\sigma_\ell<\widehat{i}\leq \sigma_{\ell+1}<\dots \sigma_r)$ above
means that we are considering
$(\sigma_0<\dots<\sigma_\ell<\sigma_{\ell+1}<\dots \sigma_r)$ but with the
additional information that $i\leq \sigma_{\ell+1}$)
from the set of faces of $\Delta^{n-1}$ to $\Delta^n_i$ induces an isomorphism
of $\omega$-categories $\Delta^{n-1}\iso \Delta^n_i$.
\ep

\bpf Both $\omega$-categories $\Delta^{n-1}$ and $\Delta^n_i$ are
freely generated by $\omega$-complexes whose atoms are clearly in
bijection. Moreover by Steiner's formulae recalled in
Section~\ref{convention}, the algebraic structure of
$\Delta^{n-1}$ and $\Delta^n_i$ is completely characterized by the
$\de^-$ and $\de^+$ operators, which are obviously preserved by
the mapping. \epf

\bp Let us denote by $I^{n+1}_j$ for $j=1,\dots,n+1$ and $n\geq 1$ the
$\omega$-subcategory of $I^{n+1}$ obtained by keeping only the words
$k_1\dots k_{n+1}$ such that $k_j=-$. Then the set map
$$\ell_1\dots \ell_n \mapsto \ell_1\dots \ell_{j-1}-\ell_j\dots \ell_n$$
from the set of faces of $I^n$ to $I^{n+1}_j$ induces an
isomorphism of $\omega$-categories $I^n \iso I^{n+1}_j$.  \ep

\bpf Same proof as for Proposition~\ref{sous-delta}. \epf

\bd One calls \textit{valid} expression of $I^{n+1}$ (resp.
$\Delta^n$) a composite of faces of $I^{n+1}$ (resp.  $\Delta^n$) of
strictly positive dimension which makes sense with respect to the
calculation rules of $I^{n+1}$.  For example $R(-0)*_0 R(0-)$ is not a
valid expression of $I^2$, whereas $R(-0)*_0 R(0+)$ is valid
(the latter being equal to $s_1R(00)$). \ed

Let $\underline{A}$ be a valid expression of $I^{n+1}$.
Suppose that $s_0 \underline{A}=R(-_{n+1})$.  Let us define by
induction on the number of faces appearing in $\underline{A}$ an
expression $\underline{A}^\bot $ using the composition laws and
variables in $\Delta^n$ (we do not know yet whether the latter
expression is valid or not in $\Delta^n$)\thinspace:
\begin{itemize}
\item if $\underline{A}$ is equal to one face $x$ of $I^{n+1}$, then this face
is necessarily in the image of $\phi_n^-$ and one can set $\underline{A}^\bot:=x^\bot$.
\item if $\underline{A}=\underline{B} *_0 \underline{C}$, then set
$\underline{A}^\bot\thinspace:= \underline{B}^\bot$.
\item if $\underline{A}=\underline{B} *_r \underline{C}$ for some $r\geq 1$,
then set $\underline{A}^\bot=\underline{B}^\bot *_{r-1} \underline{C}^\bot$.
\end{itemize}
One sees immediately by induction that whenever a face $x$ of $I^{n+1}$ such that
$s_0 x= R(-_{n+1})$ appears in $\underline{A}$, then $x^\bot$ appears in
 $\underline{A}^\bot$ because in a situation like
$\underline{A}=\underline{B} *_0 \underline{C}$, there cannot be any such face
in the expression $\underline{C}$ (since $\underline{A}$ is valid !).

If $X$ is a set of faces of $I^{n+1}$, let $X^\bot= R\left(\{y^\bot,
  y\in X\cap Im(\phi_n^-)\}\right)$.  In the case where $X$ does not
contain any element of $Im(\phi_n^-)$, then $X^\bot=\emptyset$.  If
$\underline{A}$ is a valid expression of $I^{n+1}$, then
$\underline{A}=R(\{y\hbox{ variable appearing in }\underline{A}\})$ because
$I^{n+1}$ is an $\omega$-complex.
If in addition $\underline{A}^\bot$ is a valid expression of $\Delta^n$ such
that $s_0 \underline{A}=R(-_{n+1})$,
then the two meanings of $\underline{A}^\bot$ coincide. If
$X$ and $Y$ are two sets of faces of $I^{n+1}$, then it is obvious that
$X^\bot \cup Y^\bot = (X\cup Y)^\bot$, $X^\bot \cap Y^\bot = (X\cap Y)^\bot$
and  $X^\bot \backslash Y^\bot = (X\backslash Y)^\bot$.

\bp\label{calcul2}
Let $a$ be a face of $I^{n+1}$. Then
\[\left(\de^- R(a)\right)^\bot=\de^- R(a^\bot)\] and
\[\left(\de^+ R(a)\right)^\bot=\de^+ R(a^\bot).\]
\ep

\bpf If $a$ is a $0$-dimensional face of $I^{n+1}$, then both
sides of the equalities are equal to the empty set. Let us then
suppose that $a$ is at least of dimension $1$. If $a$ does not
belong to the image of $\phi_n^-$, then both sides are again
empty. So suppose that $a\in Im(\phi_n^-)$. Then $a=k_1\dots
k_{n+1}$ where $\{k_1,\dots,k_{n+1}\}\subset \{-,0\}$. Let
$\{i_1<\dots < i_s \}=\{i\in [1,n+1],k_i=0\}$. Then
$(\de^-R(a))^\bot=R\left(\left\{k_1\dots[(-)]_{i_{2k+1}}\dots
k_{n+1},1\leq 2k+1\leq s\right\}\right)^\bot$ by definition of
$\de^-$ in $I^{n+1}$ and by definition of $(-)^\bot$ where the
notation $k_1\dots [\alpha]_i \dots k_{n+1}$ means that $k_i$ is
replaced by $\alpha$. On the other hand, $\de^- R(a^\bot)=\de^-
R((i_1-1<\dots<i_s-1))=(\de^-R(a))^\bot$. \epf

\bp\label{ortho_valid}
For any valid expression $\underline{A}$ of $I^{n+1}$ such that
$s_0 \underline{A}=R(-_{n+1})$, then $\underline{A}^\bot$ is a
valid expression of $\Delta^n$.
\ep

\bpf We are going to prove by induction on $r$ the statement
$P(r)$\thinspace: ``for any valid expression $\underline{A}$ of
$I^{n+1}$ with at most $r$ variables such that
$s_0\underline{A}=R(-_{n+1})$, the expression
$\underline{A}^\bot$ is valid in $\Delta^n$.'' If $r=1$, then
$\underline{A}=R(\{x\})$ for some face $x$ of $I^{n+1}\cap
Im(\phi_n^-)$. So  $\underline{A}^\bot=R(\{x^\bot\})$ which is
necessarily a valid expression in $\Delta^n$.
So $P(1)$ holds. Let us suppose that $P(r)$ is proved for $r\leq
r_0$ and let us consider a valid expression $\underline{A}$ of
$I^{n+1}$ with $r_0+1$ variables. Then
$\underline{A}=\underline{B}*_m \underline{C}$ with
$\underline{B}$ and $\underline{C}$ being valid expressions
having less than $r_0$ variables. By the induction hypothesis,
both $\underline{B}^\bot$ and $\underline{C}^\bot$ are valid
expressions of $\Delta^n$. If $m=0$, then by construction
$\underline{A}^\bot=\underline{B}^\bot$ and there is nothing to
prove. Otherwise by construction again,
$\underline{A}^\bot=\underline{B}^\bot *_{m-1}
\underline{C}^\bot$. Then by Steiner's formulae and by
Proposition~\ref{calcul2}, {\small
\begin{alignat*}{2}
t_{m-1}(\underline{B}^\bot) &= \left(\bigcup_{a\in \underline{B}^\bot, dim(a)\leq m-1} R(a)\right)\backslash \left(\bigcup_{b\in
\underline{B}^\bot, dim(b)=m}\left(R(b)\backslash \de^- R(b)\right)\right)&& \\
&=  \left(\bigcup_{a\in \underline{B}^\bot, dim(a)\leq m-1} R(\phi_n^-(a))^\bot\right)\backslash \left(\bigcup_{b\in
\underline{B}^\bot, dim(b)=m}\left(R(\phi_n^-(b))^\bot\backslash \de^- R(\phi_n^-(b))^\bot\right)\right)&& \\
&= \left(\bigcup_{a\in \underline{B}\cap Im(\phi_n^-), dim(a)\leq m} R(a)^\bot\right)\backslash \left(\bigcup_{b\in
\underline{B}\cap Im(\phi_n^-), dim(b)=m+1}\left(R(b)^\bot\backslash \de^- R(b)^\bot\right)\right)&&\\
&= \left(\bigcup_{a\in \underline{B}, dim(a)\leq m} R(a)^\bot\right)\backslash \left(\bigcup_{b\in
\underline{B}, dim(b)=m+1}\left(R(b)^\bot\backslash \de^- R(b)^\bot\right)\right)&&\\
&= (t_m \underline{B})^\bot &&\\
&= (s_m \underline{C})^\bot &&\\
&= s_{m-1} (\underline{C}^\bot)&&\\
\end{alignat*}
}
Consequently, $t_{m-1}\underline{B}^\bot=s_{m-1}\underline{C}^\bot$, which proves
that $\underline{A}^\bot$ is a valid expression of $\Delta^n$.
\epf

\begin{cor}
Let $\underline{A}$ and  $\underline{B}$ be two valid expressions of
$I^{n+1}$ such that $s_0 \underline{A}=s_0 \underline{B}=R(-_{n+1})$ and
such that $\underline{A}*_m \underline{B}$ exists for some $m\geq 1$.
Then $\underline{A}^\bot *_{m-1} \underline{B}^\bot$ is a valid expression
of $\Delta^n$ and
$\left(\underline{A}*_m \underline{B}\right)^\bot=\underline{A}^\bot *_{m-1} \underline{B}^\bot$.
\end{cor}

\bth\label{ortho_cube} For $n\geq 0$,
the isomorphism of $\omega$-categories
$\P^-_{R(-_{n+1})}(I^{n+1})\iso\Delta^n$ holds. \eth

\bpf The proof is threefold\thinspace: 1) one has to prove that
$\phi_n^-$ induces an $\omega$-functor $\overline{\phi}_n^-$ from
$\Delta^n$ to $\P^-(I^{n+1})$\thinspace; 2) afterwards we check
that the underlying set map of $\overline{\phi}_n^-$ is
injective\thinspace; 3) finally we prove that the image of
 $\overline{\phi}_n^-$ is the underlying set
of $\P^-_{R(-_{n+1})}(I^{n+1})$ as a whole.

\textbf{Step 1.} We are going to prove $P(n)$\thinspace: ``the set map
$\phi_n^-$ from the set of faces of $\Delta^n$ to $\P^-(I^{n+1})$ induces an
$\omega$-functor from $\Delta^n$ to $\P^-(I^{n+1})$'' by induction on
$n\geq 0$ and using Proposition~\ref{filling_simp}. The $\omega$-category
$\P^-(I^1)$ is the $\omega$-category $2_0$ generated by one $0$-morphism.
Therefore $P(0)$ holds. One has the commutative diagrams
\[
\xymatrix {
{\left\{\hbox{Faces of }\Delta_i^n\right\}} \ar@{->}[rrr]^{\phi_n^-\upharpoonright_{\Delta_i^n}}  &&& {I_{i+1}^{n+1}\subset I^{n+1}} \\
{\left\{\hbox{Faces of }\Delta^{n-1}\right\}} \fu{\iso} \ar@{->}[rrr]^{\phi_{n-1}^-} &&& {I^n} \fu{\iso}
}
\]
for any $i=0,\dots,n$.
Therefore by Proposition~\ref{filling_simp} and the induction hypothesis, one
sees that $\phi_n^-$ induces an $\omega$-functor $\overline{\phi}_n^-$ from
$\Delta^n\backslash \{(01\dots n)\}$ to $\P^-(I^{n+1})$.

It remains to prove that $s_n(R(0_{n+1}))^-= \overline{\phi}_n^-
s_{n-1} (0<\dots<n)$ and that $t_n(R(0_{n+1}))^-= \overline{\phi}_n^-
t_{n-1}(0<\dots<n)$ to complete the proof. Let us check the first
equality. By Proposition~\ref{ortho_valid}, $\left(s_n R(0_{n+1})\right)^\bot$ is a
valid expression of $\Delta^n$ in which the $(n-1)$-dimensional
faces $(0<\dots<\widehat{2i}<\dots<n)$ for $0\leq 2i\leq n$ appear. In the
$\omega$-complex $\Delta^n$, that means
that $\left(s_n R(0_{n+1})\right)^\bot$ contains not only the faces
$(0<\dots<\widehat{2i}<\dots<n)$, but also all their subfaces. Since
$(0<1<\dots<n)\notin \left(s_n R(0_{n+1})\right)^\bot$, then
necessarily $\left(s_n R(0_{n+1})\right)^\bot=s_{n-1}R(0<1<\dots<n)$.
 Because of the calculation rules in $\P^-(I^{n+1})$ described in
Theorem~\ref{calc_rules}, $\phi_n^-\left(\left(s_n R(0_{n+1})\right)^\bot\right)$ is
necessarily equal to $s_n (R(0_{n+1}))^-$. And the desired equality is then proved.

\textbf{Step 2.} The previous paragraph shows that there is a well defined $\omega$-functor $\psi_n$
from $\P^-_{R(-_{n+1})}(I^{n+1})$ to $\Delta^n$ characterized by
\begin{itemize}
\item if $x\in Im(\phi_n^-)$ (in particular $x$ is a face of $I^{n+1}$), then $\psi_n (R(x))^-= R(x^\bot)$
\item for $x,y\in \P^-_{R(-_{n+1})}(I^{n+1})$, if $x *_r y$ exists for some
$r\geq 1$, then $\psi_n(x *_r y)=\psi_n(x) *_{r-1} \psi_n(y)$.
\end{itemize}
Moreover $\psi_n \circ \overline{\phi}_n^-=\id$ and therefore $\overline{\phi}_n^-$ is injective.

\textbf{Step 3.} This is an immediate consequence of Theorem~\ref{calc_rules} and of the construction
of $\overline{\phi}_n^-$.
\epf

\begin{cor}\label{ortho_complet}
For $n\geq 0$,
\[\P^-(I^{n+1})=\bigoplus_{1\leq p\leq n+1} (\Delta^{p-1})^{\oplus C^p_{n+1}}\]
where $\oplus$ is the direct sum in the category  of $\omega$-categories (which
corresponds for the underlying sets to the disjoint union).
\end{cor}

\bpf By Proposition~\ref{grading}, the $\omega$-category $\P^-(I^{n+1})$ is graded
by the vertices of $I^{n+1}$ (however $G^{+_{n+1},-}\P^-(I^{n+1})=\emptyset$).
So this is a consequence of Theorem~\ref{ortho_cube}.
\epf

The isomorphism of Theorem~\ref{ortho_cube} is actually more than only
an isomorphism of $\omega$-categories. Indeed, let
$\underline{\Delta}$ be the unique small category such that a presheaf
over $\underline{\Delta}$ is exactly a simplicial set
\cite{May,Weibel}. The category $\underline{\Delta}$ has for objects
the finite ordered sets $[n]=\{0<1<\dots<n\}$ for all integers $n\geq
0$ and has for morphisms the nondecreasing monotone functions. One is
used to distinguish in this category the morphisms
$\epsilon_i:[n-1]\rightarrow [n]$ and $\eta_i:[n+1]\rightarrow [n]$
defined as follows for each $n$ and $i=0,\dots,n$\thinspace:
\[
\epsilon_i(j)=\left\{\begin{array}{c}j \hfill\hbox{ if }j<i\\j+1\hfill\hbox{ if }j\geq i\end{array}\right\}
,\ \ \
\eta_i(j)=\left\{\begin{array}{c}j \hfill\hbox{ if }j\leq i\\j-1\hfill\hbox{ if }j> i\end{array}\right\}
\]

It is well-known (\cite{oriental}) that the map $[n]\mapsto \Delta^n$ induces a
functor $\Delta^*$ from
$\underline{\Delta}$ to $\omega Cat$ by setting
$\epsilon_i \mapsto \Delta^{\epsilon_i}$ and $\eta_i \mapsto \Delta^{\eta_i}$ where
\begin{itemize}
\item for any face $(\sigma_0<\dots <\sigma_s)$ of $\Delta^{n-1}$,
$\Delta^{\epsilon_i}(\sigma_0<\dots <\sigma_s)$ is the
only face of $\Delta^{n}$ having $\epsilon_i\{\sigma_0,\dots,\sigma_s\}$ as set of vertices\thinspace;
\item for any face $(\sigma_0<\dots <\sigma_r)$ of $\Delta^{n+1}$,
$\Delta^{\eta_i}(\sigma_0<\dots <\sigma_r)$ is the
only face of $\Delta^{n}$ having $\eta_i\{\sigma_0,\dots,\sigma_r\}$ as set of vertices.
\end{itemize}

For a face $\sigma=\{\sigma_1,\dots,\sigma_s\}$ of $\Delta^n$ with $\sigma_1<\dots<\sigma_s$,
let $\phi_n^-(\sigma):=\phi_n^-(\sigma_1<\dots<\sigma_s)$. Then one has

\bp\label{cube-simplexe} Let $\delta_i^-:I^{n}\rightarrow
I^{n+1}$ be the $\omega$-functor corresponding by Yoneda's lemma
to the face map $\de_i^-:\omega Cat(I^{n+1},-)\rightarrow \omega
Cat(I^{n},-)$ and let $\gamma_i^-:I^{n+2}\rightarrow I^{n+1}$ be
the $\omega$-functor corresponding by Yoneda's lemma to the
degeneracy maps $\Gamma_i^-:\omega Cat(I^{n+1},-)\linebreak[4]
\rightarrow \omega Cat(I^{n+2},-)$ for $1\leq i\leq n+1$. Then
the following diagrams are commutative\thinspace:
\[
\xymatrix{
{\Delta^{n-1}}\fr{\Delta^{\epsilon_i}} \fd{\phi_{n-1}^-} & {\Delta^{n}}\fd{\phi_{n}^-}\\
{I^n} \fr{\delta_{i+1}^-} & {I^{n+1}}
}\ \ \
\xymatrix{
{\Delta^{n+1}}\fr{\Delta^{\eta_i}} \fd{\phi_{n+1}^-} & {\Delta^{n}}\fd{\phi_{n}^-}\\
{I^{n+2}} \fr{\gamma_{i+1}^-} & {I^{n+1}}
}
\]
\ep

\bpf This is a remake of the proof that the natural morphism
$h^-$ from the globular to the branching semi-cubical nerve preserves
the simplicial structure of both sides (see \cite{sglob}). \epf

Remember that the geometric explanation of Proposition~\ref{cube-simplexe} is again that
close to a corner, the intersection of an $n$-cube  by an hyperplane is exactly the
$(n-1)$-simplex.

\begin{cor} For $0\leq i\leq n$, the mapping
\[k_1\dots k_{n+1}\mapsto \delta_{i+1}^-(k_1\dots k_{n+1})= k_1\dots k_i [-]_{i+1} k_{i+1}\dots k_{n+1}\]
induces an $\omega$-functor from $\P^-_{R(-_{n})}(I^{n})$ to
$\P^-_{R(-_{n+1})}(I^{n+1})$. And the mapping
\[k_1\dots k_{n+2}\mapsto \gamma_{i+1}^-(k_1\dots k_{n+2})= k_1\dots k_i max(k_{i+1},k_{i+2})k_{i+3}\dots k_{n+2}\]
induces an $\omega$-functor from  $\P^-_{R(-_{n+2})}(I^{n+2})$ to
$\P^-_{R(-_{n+1})}(I^{n+1})$. This way, the mapping $[n]\mapsto \P^-_{R(-_{n+1})}(I^{n+1})$ induces
a functor from $\underline{\Delta}$ to $\omega Cat$ which is isomorphic to the functor $\Delta^*$.
In other terms, the family of maps $(\phi_n)_{n\geq 0}$ induces an isomorphism of
functors from $\underline{\Delta}$ to $\omega Cat$
\[\Delta^* \iso \P^-_{R(-_{*+1})}(I^{*+1}).\]
\end{cor}

\section{Simplicial cut}\label{genera1}

The branching semi-globular nerve will not match the notion of
cut as presented in \cite{sglob}. A slight modification of that
definition is therefore necessary. This is precisely the purpose of
this section. The reader does not need to know the previous definition
of a cut.

\bd A cut is a triple $(\mathcal{F},\underline{P},\ev)$ where
$\underline{P}$ is a functor from $\omega Cat_1$ to $\omega Cat$,
$\mathcal{F}$ a functor from $\omega Cat_1$ to $Sets^{\Delta^{op}}_+$,
and $\ev=(\ev_n)_{n\geq 0}$ a family of natural transformations
$\ev_n:\mathcal{F}_n\longrightarrow tr^n \underline{P}$ where
$\mathcal{F}_n$ is the set
of $n$-simplexes of $\mathcal{F}$.  A morphism of cuts from
$(\mathcal{F},\underline{P},\ev)$ to $(\mathcal{G},\underline{Q},\ev)$
is a pair $(\phi,\psi)$ such that $\phi$ is a natural transformation
of functors from $\mathcal{F}$ to $\mathcal{G}$ and $\psi$ a natural
transformation of functors from $\underline{P}$ to $\underline{Q}$
which makes the following diagram commutative for any $n\geq 0$\thinspace:
\[\xymatrix{{\mathcal{F}_n} \fd{\phi_n} \fr{\ev_n}& {tr^n \underline{P}}\fd{\psi_n}\\
{\mathcal{G}_n} \fr{\ev_n} & {tr^n \underline{Q}}}\]
\ed

If $(\mathcal{F},\ev)$ is a cut in the sense of \cite{sglob}, then
$(\mathcal{F},\P,\ev)$ is a cut in the sense of the above
definition. One can omit $\ev$ and simply denote the cut
$(\mathcal{F},\underline{P},\ev)$ by
$(\mathcal{F},\underline{P})$.

For such a cut, one can define the associated homology theory as
in \cite{sglob}. If $(\mathcal{F},\underline{P})$ is a cut, let
$C_{n+1}^{(\mathcal{F},\underline{P})}(\C):=C_n(\mathcal{F}(\C))$
and let $H_{n+1}^{(\mathcal{F},\underline{P})}$ be the
corresponding homology theory for $n\geq -1$.

To illustrate the definition, let us recall now the definition of the
branching semi-cubical cut \cite{Gau,Coin}. Let $\C$ be a non-contracting
$\omega$-category. We set
{\small
\[\omega Cat(I^{n+1},\C)^{-}:=\left\{x\in
\omega Cat(I^{n+1},\C),d_0^{-}(u)={-}_{n+1} \hbox{ and }dim(u)=1 \Longrightarrow dim(x(u))=1\right\}\]}
where ${-}_{n+1}$ is the initial state
 of $I^{n+1}$. For all
$(i,n)$ such that  $0\leq i\leq n$, the face maps $\de_i$ from
$\omega Cat(I^{n+1},\C)^{-}$ to $\omega Cat(I^{n},\C)^{-}$
are the arrows $\de^{-}_{i+1}$ defined by
\[\de^{-}_{i+1}(x)(k_1\dots k_{n+1})=x(k_1\dots [{-}]_{i+1}\dots k_{n+1})\]
and the degeneracy maps $\epsilon_i$ from $\omega Cat(I^{n},\C)^{-}$
to $\omega Cat(I^{n+1},\C)^{-}$
are the arrows $\Gamma^{-}_{i+1}$ defined by setting
\[\Gamma_i^-(x)(k_1\dots k_n):=x(k_1\dots \max(k_i,k_{i+1})\dots k_n)\]
with the order $-<0<+$.

Let $\C$ be a non-contracting $\omega$-category. The $\N$-graded set
$\mathcal{N}^{-}_*(\C)=\omega Cat(I^{*+1},\C)^-$ together with the
convention $\mathcal{N}^{-}_{-1}(\C)=\C_0$, endowed with the maps
$\de_i$ and $\epsilon_i$ above defined with moreover $\de_{-1}=s_0$
and with $\ev(x)=x(0_n)$ for $x\in \omega Cat(I^{n},\C)$ is an
augmented simplicial set and $(\mathcal{N}^-,\P)$ becomes a simplicial
cut. It is called the \textit{branching simplicial cut} associated to
$\C$. Set $H_{n+1}^{{-}}(\C):=H_n(\mathcal{N}^{{-}}(\C))$ for $n\geq -1$.
This homology theory is called the \textit{branching semi-cubical
homology}.

\section{Thin elements in simplicial cuts and reduced homology}\label{genera2}

Let $(\mathcal{F},\underline{P})$ be a simplicial cut.  Let
$M_n^\mathcal{F}:\omega Cat_1\longrightarrow Ab$ be the functor
defined as follows\thinspace: the group
$M_{n}^{(\mathcal{F},\underline{P})}(\C)$ is the subgroup generated by
the elements $x\in \mathcal{F}_{n-1}(\C)$ such that $\ev(x)\in
tr^{n-2}\underline{P}\C$ for $n\geq 2$ and with the convention
$M_0^{(\mathcal{F},\underline{P})}(\C)=M_1^{(\mathcal{F},\underline{P})}(\C)=0$
and the definition of $M_n^\mathcal{F}$ is obvious on non-contracting
$\omega$-functors.

\bd The elements of $M_{*}^{(\mathcal{F},\underline{P})}(\C)$ are called
\textit{thin}. \ed

To illustrate the definition above, let us consider again the case of
the branching semi-cubical cut. A thin element of
$\mathcal{N}_{n-1}^-(\C)$ is nothing else but an $\omega$-functor
$f\in\omega Cat(I^{n},\C)^{-}$ such that $f(0_{n})\in tr^{n-2}\P\C=tr^{n-1}\C$,
that is $f(0_{n})$ is of dimension at most $n-1$. So $f$ corresponds intuitively
to a $n$-cube without volume.

Let us come back now to the general situation. Let
$\CR^{(\mathcal{F},\underline{P})}_n:\omega Cat_1\longrightarrow
\comp$ be the functor defined by
$$\CR^{(\mathcal{F},\underline{P})}_n:=C_n^{(\mathcal{F},\underline{P})}/(M_n^{(\mathcal{F},\underline{P})}+\de
M_{n+1}^{(\mathcal{F},\underline{P})})$$ where $\comp$ is the category
of chain complexes of abelian groups and endowed with the differential
map $\de$.

\bd This chain complex is called the \textit{reduced complex}
associated to the cut ${(\mathcal{F},\underline{P})}$ and the
corresponding homology is denoted by
$\HR^{(\mathcal{F},\underline{P})}_*$ and is called the
\textit{reduced homology} associated to
${(\mathcal{F},\underline{P})}$. \ed

A morphism of cuts from ${(\mathcal{F},\underline{P})}$ to
${(\mathcal{G},\underline{Q})}$ yields natural morphisms from
$H_*^{(\mathcal{F},\underline{P})}$ to
$H_*^{(\mathcal{G},\underline{Q})}$ and from
$\HR_*^{(\mathcal{F},\underline{P})}$ to
$\HR_*^{(\mathcal{G},\underline{Q})}$. There is also a canonical
natural transformation $R^{(\mathcal{F},\underline{P})}$ from
$H_*^{(\mathcal{F},\underline{P})}$ to
$\HR_*^{(\mathcal{F},\underline{P})}$, functorial with respect to
${(\mathcal{F},\underline{P})}$, that makes the following diagram
commutative\thinspace:
\[
\xymatrix{{H_*^{(\mathcal{F},\underline{P})}}\fr{R^{(\mathcal{F},\underline{P})}}\fd{} & {\HR_*^{(\mathcal{F},\underline{P})}}\fd{}\\ {H_*^{(\mathcal{G},\underline{Q})}}\fr{R^(\mathcal{G},\underline{Q})}& {\HR_*^{(\mathcal{G},\underline{Q})}}}
\]

\begin{nota} The reduced branching semi-cubical complex is denoted by
$\CR^-$ and the reduced branching semi-cubical homology theory is
denoted by $\HR^{-}$. \end{nota}

\section{Regular cut}\label{genera3}

The next definition presents the notion of \textit{regular cut}.
It makes the construction of folding operators possible.

\bd\label{def_regular}
A cut ${(\mathcal{F},\underline{P})}$ is \textit{regular} if and only if it satisfies the following
properties\thinspace:
\begin{enumerate}
\item\label{regular0} for any $\omega$-category $\C$, the set $\mathcal{F}_{-1}(\C)$ only depends on $tr^0\C=\C_0$\thinspace:
i.e. for any $\omega$-categories $\C$ and $\D$, $\C_0=\D_0$ implies $\mathcal{F}_{-1}(\C)=\mathcal{F}_{-1}(\D)$.
\item\label{regular1} $\mathcal{F}_0:=tr^0 \underline{P}$.
\item \label{regular1.5} $\ev \circ \epsilon_i=\ev$.
\item\label{regular2} for any natural transformation of functors $\mu $ from $\mathcal{F}_{n-1}$ to
 $\mathcal{F}_{n}$ with $n\geq 1$, and for any natural map $\square$ from
$tr^{n-1}\underline{P}$ to $\mathcal{F}_{n-1}$ such that
$\ev\circ \square=\id_{tr^{n-1}\underline{P}}$, there exists one and only one
natural transformation $\mu .\square$ from $tr^n\underline{P}$ to $\mathcal{F}_{n}$ such
that the following diagram commutes
\[\xymatrix{\ar@/^20pt/[rr]^{\id_{tr^n\underline{P}}}{tr^n \underline{P}}\fr{\mu .\square}& {\mathcal{F}_n}\fr{\ev_n}&{tr^n \underline{P}}\\
\ar@/_20pt/[rr]_{\id_{tr^{n-1}\underline{P}}}{tr^{n-1}\underline{P}}\fu{i_{n}}\fr{\square}&{\mathcal{F}_{n-1}}\fr{\ev_{n-1}}\fu\mu &{tr^{n-1}\underline{P}}\fu{i_{n}}}\]
where $i_{n}$ is the canonical inclusion functor from
$tr^{n-1}\underline{P}$ to $tr^n \underline{P}$.
\item\label{regular3}  let $\square_1^{(\mathcal{F},\underline{P})}:=\id_{\mathcal{F}_0}$ and
$\square_n^{(\mathcal{F},\underline{P})}:=\epsilon_{n-2}.\dots\epsilon_0
. \square_1^{(\mathcal{F},\underline{P})}$ be natural
transformations from $tr^{n-1}\underline{P}$ to
$\mathcal{F}_{n-1}$ for $n\geq 2$\thinspace; then the natural
transformations $\de_i\square_n^{(\mathcal{F},\underline{P})}$
for $0\leq i \leq n-1$ from $tr^{n-1}\underline{P}$ to
$\mathcal{F}_{n-2}$ satisfy the following properties
\begin{enumerate}
\item\label{regular3.01} $\left\{\ev \de_{n-2}\square_n^{(\mathcal{F},\underline{P})}, \ev \de_{n-1}\square_n^{(\mathcal{F},\underline{P})}\right\}=\left\{s_{n-1},t_{n-1}\right\}$.
\item\label{regular3.02} if for some $\omega$-category $\C$ and some $u\in \C_n$, $\ev\de_i\square_n^{(\mathcal{F},\underline{P})}(u)=d_p^\alpha u$
for some $p\leq n$ and for some $\alpha\in\{-,+\}$, then $\de_i\square_n^{(\mathcal{F},\underline{P})}=\de_i\square_n^{(\mathcal{F},\underline{P})} d_p^\alpha $.
\end{enumerate}
\item\label{regular3.1} let  $\Phi_n^{(\mathcal{F},\underline{P})}:=\square_n^{(\mathcal{F},\underline{P})}\circ \ev$ be a natural
transformation from $\mathcal{F}_{n-1}$ to itself\thinspace; then $\Phi_n^{(\mathcal{F},\underline{P})}$
induces the identity natural transformation  on $\CR_{n}^{(\mathcal{F},\underline{P})}$.
\item\label{regular4} if $x$, $y$ and $z$ are three elements of $\mathcal{F}_n(\C)$,
and if $\ev(x)*_p \ev(y)=\ev(z)$ for some $1\leq p\leq n$, then
$x+y=z$ in $\CR_{n+1}^{{(\mathcal{F},\underline{P})}}(\C)$ and in a functorial way.
\end{enumerate}
If ${(\mathcal{F},\underline{P})}$ is a regular cut, then the
natural transformation $\Phi_n^{(\mathcal{F},\underline{P})}$ is
called the $n$-dimension\-al folding operator of the cut
${(\mathcal{F},\underline{P})}$. By convention, one sets
$\square_0^{(\mathcal{F},\underline{P})}=\id_{\mathcal{F}_{-1}}$
and $\Phi_0^{(\mathcal{F},\underline{P})}=\id_{\mathcal{F}_{-1}}$.
There are no ambiguities to set
$\Phi^{(\mathcal{F},\underline{P})}(x):=\Phi^{(\mathcal{F},\underline{P})}_{n+1}(x)$
for $x\in \mathcal{F}_n(\C)$ for some $\omega$-category $\C$. So
$\Phi^{(\mathcal{F},\underline{P})}$ defines a natural
transformation, and even a morphism of cuts, from
${(\mathcal{F},\underline{P})}$ to itself. However beware of the
fact that there is really an ambiguity in the notation
$\square^{(\mathcal{F},\underline{P})}$\thinspace: so the latter
will not be used. \ed

Now here are some trivial remarks about regular cuts\thinspace:
\begin{itemize}
\item Let $f$ be a natural set map from $tr^0 \underline{P}\C=\C_1$ to itself. Let $2_1$ be the
$\omega$-category generated by one $1$-morphism $A$. Then necessarily
$f(A)=A$ and therefore $f=\id$. So the above axioms imply that
$\ev_0=\id$.

\item The map $\Phi_n^{(\mathcal{F},\underline{P})}$ induces the identity natural transformation  on $\HR^{(\mathcal{F},\underline{P})}_{n}$.
\item For any $n\geq 1$, there exists non-thin elements $x$ in $\mathcal{F}_{n-1}(\C)$
as soon as $\C_n\neq \emptyset$. Indeed, if $u\in \C_n$,  $\ev\ \square_n^{(\mathcal{F},\underline{P})}(u)=u$,
therefore $\square_n^{(\mathcal{F},\underline{P})}(u)$ is a non-thin element of $\mathcal{F}_{n-1}(\C)$.
\end{itemize}

We end this section by some general facts about regular cuts.

\bp\label{Phi_functoriel} Let $f$ be a morphism of cuts from ${(\mathcal{F},\underline{P})}$ to ${(\mathcal{G},\underline{Q})}$. Suppose that ${(\mathcal{F},\underline{P})}$ and
${(\mathcal{G},\underline{Q})}$ are regular. Then the equality
$\Phi^{(\mathcal{G},\underline{Q})} \circ f=f\circ \Phi^{(\mathcal{F},\underline{P})}$ holds as natural
transformation from ${(\mathcal{F},\underline{P})}$ to ${(\mathcal{G},\underline{Q})}$. In other terms, the following diagram
is commutative\thinspace:
\[ \xymatrix{{(\mathcal{F},\underline{P})} \fr{f} \fd{\Phi^{(\mathcal{F},\underline{P})}}& {(\mathcal{G},\underline{Q})} \ar@{->}[d]^{\Phi^{(\mathcal{G},\underline{Q})}} \\
{(\mathcal{F},\underline{P})} \fr{f} & {(\mathcal{G},\underline{Q})}} \]
\ep

\bpf See \cite{sglob}.
\epf

\bp\label{canonical_position} If $u$ is a $(n+1)$-morphism of $\C$ with $n\geq 1$,
then $\square_{n+1}^{(\mathcal{F},\underline{P})} u$ is an homotopy within the
simplicial set $\mathcal{F}(\C)$ between $\square_n^{(\mathcal{F},\underline{P})} s_nu$
and $\square_n^{(\mathcal{F},\underline{P})} t_nu$. \ep

\bpf See \cite{sglob}.
\epf

\begin{cor}\label{diff_reduite}
If  $x\in \CR^{(\mathcal{F},\underline{P})}_{n+1}(\C)$, then
$\de x=\de\square_{n+1}^{(\mathcal{F},\underline{P})}x=\square_{n}^{(\mathcal{F},\underline{P})} s_nx-\square_{n}^{(\mathcal{F},\underline{P})}t_nx$ in
$\CR^{(\mathcal{F},\underline{P})}_n(\C)$.  In other terms, the differential map from
$\CR^{(\mathcal{F},\underline{P})}_{n+1}(\C)$ to $\CR^{(\mathcal{F},\underline{P})}_{n}(\C)$ with $n\geq 1$ is induced by the
map $s_n-t_n$.
\end{cor}

Now let us come back to the particular example of the branching semi-cubical
cut.

\bth \cite{Coin}
The branching semi-cubical cut is regular. \eth

\begin{nota} The branching semi-cubical folding operators are denoted by
$\Phi^{-}$ and $\square_n^{-}$. \end{nota}

\bd\cite{Coin}\label{formalcoin}  Set
\begin{itemize}
\item $\CF^{-}_0(\C):=\Z\C_0$
\item $\CF^{-}_1(\C):=\Z\C_1$
\item $\CF^{-}_n(\C)= \Z\C_n/\{x*_0y=x,x*_1y=x+y, \dots , x
*_{n-1}y=x+y\hbox{ mod }\Z tr^{n-1}\C\}$ for $n\geq 2$
\end{itemize}
with the differential map $s_{n-1}-t_{n-1}$ from $\CF^{-}_n(\C)$ to
$\CF^{-}_{n-1}(\C)$ for $n\geq 2$ and $s_0$ from
$\CF^{-}_1(\C)$ to $\CF^{-}_0(\C)$. This chain complex is called
the formal negative corner complex. The associated homology is denoted by
$\HF^{-}(\C)$ and is called the formal negative corner homology of $\C$.
The map $\CF_*^-$ (resp. $\HF_*^-$) induces a functor from $\omega Cat_1$
to the category of chain complexes of groups  $\comp$ (resp. to the category of abelian
groups $Ab$).
\ed

Since the branching semi-cubical cut is regular, then
$\square^-_n(x*_py)=\square^-_n(x)+\square^-_n(y)$ in
$\CR_n^-(\C)$ for any non-contracting $\omega$-category $\C$, for
any $p\geq 1$ and any $x,y\in \C_n$ as soon as $x*_py$ exists.
Moreover we have with the same notation
$\square^-_n(x*_0y)=\square^-_n(x)$ as soon as $x*_0 y$ exists
\cite{Coin}. So the folding operators $\square^-_n$ induce a
natural morphism of chain complexes from $\CF^-$ to $\CR^-$.

\section{The branching semi-globular nerve}\label{definition}

Using the functor $\Delta^*$ from $\underline{\Delta}$ to $\omega
Cat$, one obtains the well-known simplicial nerve of
$\omega$-category $\mathcal{N}_*(\C):=\omega Cat(\Delta^*,\C)$
introduced for the first time in \cite{oriental}. Now we have the
necessary tools in hand  to define the branching semi-globular cut
of a non-contracting $\omega$-category.

\bd Let $\C$ be a non-contracting $\omega$-category. Then set
$$\mathcal{N}^{gl^-}_n(\C)=\omega Cat(\Delta^n,\P^-\C)$$
and $\mathcal{N}^{gl^-}_{-1}(\C)=\C_0$ with $\de_{-1}(x)=s_0 x$.
Then $\mathcal{N}^{gl^-}$ induces a functor from $\omega Cat_1$ to
$Sets^{\Delta^{op}}_+$. The triple $(\mathcal{N}^{gl^-},\P^-,\ev)$
with $\ev(x)=x(0<1<\dots<n)$ is called the branching semi-globular cut.
\ed

\bp If $\P\C$ is a strict globular $\omega$-groupoid, then the
branching semi-globu\-lar cut of $\C$ is Kan. \ep

\bpf This is an immediate corollary of Theorem~\ref{groupoid_aussi}.  \epf

\bth The branching semi-globular cut is regular. \eth

\bpf Let $\C$ be a non-contracting $\omega$-category.  Let $\D$ be the
unique $\omega$-category such that $\P\D=\P^-\C$ and with
$s_0(\P\D)=\{\alpha\}$, $t_0(\P\D)=\{\beta\}$ and $\alpha\neq \beta$.
Then $\mathcal{N}_n^{gl}(\D)=\mathcal{N}_n^{gl^-}(\C)$ for any $n\geq
0$.  So the regularity of the branching semi-globular nerve comes from that of the
globular nerve \cite{sglob}.  \epf

The reduced branching semi-globular homology theory is denoted by $\HR^{gl^-}$, the
branching semi-globular folding operators by $\Phi^{gl^-}$ and $\square_n^{gl^-}$.

\bth\label{relation_left_glob} Let $\C$ be a non-contracting $\omega$-category. Then for any $n\geq 0$,
$\CR_n^{gl^-}(\C)$ is generated by the $\Phi^{gl^-}(h^-(x))$ for $x$ running
over $\C_n$. Moreover, one has
\begin{itemize}
\item if $x,y\in \C_n$ and if $x*_0y$ exists, then $\Phi^{gl^-}(h^-(x*_0 y))=\Phi^{gl^-}(h^-(x))$
\item if $x,y\in \C_n$ and if $x*_ry$ exists for $r\geq 1$,
then $\Phi^{gl^-}(h^-(x*_r y))=\Phi^{gl^-}(h^-(x))+ \Phi^{gl^-}(h^-(y))$ in $\CR_n^{gl^-}(\C)$.
\end{itemize}
\eth

\bpf The group $\CR_n^{gl^-}(\C)$ is generated by the $\Phi^{gl^-}(X)$
for $X\in\P^-\C$.  By Theorem~\ref{calc_rules},
$X=\Psi(h^-(x_1),\dots,h^-(x_s))$ where $\Psi$ is an expression using
only the composition laws $*_h$ for $h\geq 1$ and where
$x_1,\dots,x_s\in \P\C$. By regularity of the branching semi-globular nerve, one deduces that
$\CR_n^{gl^-}(\C)$ is generated by the $\Phi^{gl^-}(h^-(x))$ for $x$ running
over $\C_n$. The remaining part of the statement is clear. \epf

So, as the family of operators $\square_n^-$ does for the semi-cubical
case, the family of operators $(\square_n^{gl^-})_{n\geq 0}$ induces a
natural morphism of chain complexes from $\CF^-_*$ to $\CR^{gl^-}_*$
by Theorem~\ref{relation_left_glob}.

We close this section by the statement of the ``thin elements''
conjecture for the branching semi-globular nerve\thinspace:

\begin{conj}
Let $\C$ be a non-contracting $\omega$-category which is the free
$\omega$-category freely generated by a precubical set. Consider a
linear combination $\sum_{i\in I}\lambda_i x_i$ of elements of
$M_n^{gl^-}(\C)$ with $\forall i\in I, \lambda_i\in \Z$.  Then this
linear combination is a cycle if and only if it is a boundary for the
simplicial differential map. \end{conj}

\begin{conj}\cite{Coin} (The thin elements conjecture)
Let $K$ be a precubical set. Then the chain complex morphisms
$C^{gl^\pm}(F(K))\rightarrow
\CR^{gl^\pm}(F(K))$ induce isomorphisms in homology. \end{conj}

\begin{conj} (The extended  thin elements conjecture)
Let $K$ be a precubical set.  Then the chain complex morphisms
$C^{gl^\pm}(\widehat{F(K)})\rightarrow \CR^{gl^\pm}(\widehat{F(K)})$ induce
isomorphisms in homology. \end{conj}

\section{Comparison with the branching semi-cubical nerve}\label{comparison}

\begin{figure}
\[
\xymatrix{
&&& \mathcal{N}^{gl} \ar@{->}[drrr]^{h^+} \ar@{->}[dlll]_{h^-} \ar@{->}[ddrrr]_{\mathcal{N}(\P\C\rightarrow \P^+\C)}  \ar@{->}[ddlll]^{\mathcal{N}(\P\C\rightarrow \P^-\C)} &&& \\
\mathcal{N}^{-} \fd{f\mapsto f^-}  &&&&&& \mathcal{N}^{+} \ar@{->}[d]^{f\mapsto f^+} \\
\mathcal{N}^{gl^-}   &&&&&& \mathcal{N}^{gl^+}
}
\]
\caption{Recapitulation of all simplicial constructions}
\label{recap}
\end{figure}
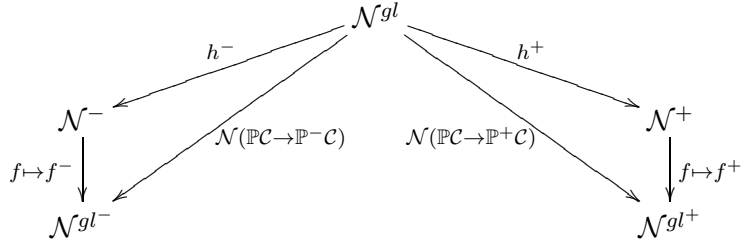

\bth\label{kani}
There exists one and only one morphism of cuts from the branching
semi-cubical cut to the branching semi-globular cut such that the
underlying natural transformation from $\P$ to $\P^-$ is the canonical
morphism $\P\rightarrow \P^-$ appearing in the definition of
$\P^-$. \eth

\bpf Let $f_1$ and $f_2$ be two morphisms of cuts from the branching
semi-cubical nerve to the branching semi-globular nerve. For any
non-contracting $\omega$-category $\C$, $(f_1)_{-1}$ and $(f_2)_{-1}$
induce natural set maps from $\C_0$ to $\C_0$. Since the only natural
transformation from the identity functor of the category of sets to
itself is the identity transformation (to see that, consider the case
of a singleton), then $(f_1)_{-1}$ and $(f_2)_{-1}$ are both equal to
the identity of $\C_0$. If $\C$ is the free $\omega$-category $2_1(A)$
generated by one $1$-morphism $A$, then $\C$ is non-contracting. In
that case, $\mathcal{N}^-_{0}(\C)=\{A\}= \mathcal{N}^{gl^-}_{0}(\C)$.
So in that case, $(f_1)_0(A)=(f_2)_0(A)$.  For any non-contracting
$\omega$-category $\C$ and any $x\in \C_1$, there exists a unique
$\omega$-functor $\underline{x}$ from $2_1(A)$ to $\C$ such that
$\underline{x}(A)=x$. So by naturality, $(f_1)_0(x)=(f_2)_0(x)$, and
therefore $(f_1)_0$ and $(f_2)_0$ coincide everywhere. Notice that in
general $\mathcal{N}^{gl^-}_{0}(\C)\neq \C_1$ so the reasoning of the
$0$-th dimension does not apply to dimension $1$. Now suppose that we
have proved that $(f_1)_n=(f_2)_n$ for $n\leq n_0$ and $n_0\geq 1$.
Then $(f_1)_{n_0+1}$ and $(f_2)_{n_0+1}$ are two $\omega$-functors
from $\Delta^{n_0+1}$ to $\P^-\C$ which coincide on the
$n_0$-dimensional faces of $\Delta^{n_0+1}$ and such that
$(f_1)_{n_0+1}(0<\dots <n_0+1)=(f_2)_{n_0+1}(0<\dots <n_0+1)$.  Then
$(f_1)_{n_0+1}$ and $(f_2)_{n_0+1}$ induce the same labeling of the
faces of $\Delta^{n_0+1}$, therefore by freeness of $\Delta^{n_0+1}$,
$(f_1)_{n_0+1}=(f_2)_{n_0+1}$.

Now let us prove the existence of this natural transformation.
Let $f\in \omega Cat(I^{n+1},\C)^-$. Then for any $x\in
G^{{R(-_{n+1})},-}\P I^{n+1}$, the morphism $f(x)$ cannot be
$0$-dimensional otherwise $s_1f(x)=f(s_1x)$ would be so as well.
So the restriction of $f$ to ${G^{{R(-_{n+1})},-}\P
I^{n+1}}\subset \P I^{n+1}$ gives rise to an element of
$$\mathcal{N}^1_n(\C):=\omega Cat\left({G^{{R(-_{n+1})},-}\P
I^{n+1}},\P\C\right)$$ and therefore to elements of
$$\mathcal{N}^2_n(\C):=\omega Cat\left({G^{{R(-_{n+1})},-}\F \left(\U \P I^{n+1}\right),{\F\left(\U \P\C\right)}}\right)$$
and of
$$\mathcal{N}^3_n(\C):=\omega Cat\left({G^{{R(-_{n+1})},-}\F \left(\U \P I^{n+1}/\Rm\right),{\F\left(\U \P\C/\Rm\right)}}\right).$$

\begin{figure}
{\small
\[
\xymatrix{
       & {\F\left(\U \P\C\right)} \ar@{->}[rr]\ar@{->}'[d][dd]
         & & {\F \left(\U \P\C/\Rm\right)} \ar@{->}[dd]
      \\
       {G^{{R(-_{n+1})},-}\F \left(\U \P I^{n+1}\right)} \ar@{->}[ur]\ar@{->}[rr]\ar@{->}[dd]
       & & {G^{{R(-_{n+1})},-}\F \left(\U \P I^{n+1}/\Rm\right)} \ar@{->}[ur]\ar@{->}[dd]
      \\
       & {\P\C} \ar@{->}'[r][rr]
         & & \cocartesien {\P^-\C}
      \\
       {G^{{R(-_{n+1})},-}\P I^{n+1}} \ar@{->}[rr]\ar@{->}[ur]
       & & \cocartesien {\P^-_{R(-_{n+1})} I^{n+1}} \ar@{-->}[ur]_{f^-}
      }
\]}
\caption{The canonical morphism from the branching semi-cubical to the branching semi-globular nerve}
\label{fmoins}
\end{figure}

The $\omega$-functors $\delta_i^-:I^n\rightarrow I^{n+1}$ and
$\gamma_i^-:I^{n+2}\rightarrow I^{n+1}$ for $1\leq i\leq n+1$ are all
non-contracting. Since the natural maps $\mathcal{N}^-_n\rightarrow
\mathcal{N}^i_n$ for $i=1,2,3$ arise from restrictions, then one
obtains three natural morphisms of simplicial sets
$\mathcal{N}^-\rightarrow \mathcal{N}^i$ which yield a cone based on
the diagram
\[
\xymatrix{
& {\mathcal{N}^3} \fd{}\\
{\mathcal{N}^1} \fr{} & {\mathcal{N}^2}
}
\]
Therefore one obtains a natural transformation
\[\mathcal{N}^-\longrightarrow \varprojlim \mathcal{N}^i \iso \omega Cat\left({\P^-_{R(-_{*+1})} I^{*+1}},-\right)
\iso \mathcal{N}^{gl^-}\]
For $f\in \mathcal{N}^-(\C)$, the corresponding element $f^-\in \mathcal{N}^{gl^-}$
is represented in Figure~\ref{fmoins}.
\epf

Let $\C$ be a non-contracting $\omega$-category. The natural
transformation $\mathcal{N}^-\rightarrow\mathcal{N}^{gl^-}$ yields a
natural morphism of chain complexes $\CR^-_*(\C)\rightarrow
\CR^{gl^-}_*(\C)$ by Section~\ref{genera2}.

\bp The following diagram is commutative and the three maps are surjective\thinspace:
\[
\xymatrix{
& {\CF^-_*} \ar@{->}[dr]^{\square_*^{gl^-}} \ar@{->}[dl]_{\square_*^{-}}& \\
{\CR^-_*} \ar@{->}[rr] && {\CR^{gl^-}_*}}
\]
\ep

\bpf Let $\C$ be a non-contracting $\omega$-category. The only thing one has
to prove is the surjectivity of $\CR^-_*(\C)\rightarrow
\CR^{gl^-}_*(\C)$.  There is nothing to prove for $n=0$ and $n=1$. So
let us suppose that $n\geq 2$. The group $\CR^{gl^-}_n(\C)$ is
generated by the elements of the form $\square_n^{gl^-}(x)$ where
$x\in (\P^-\C)_{n-1}$. The canonical map $\P\C\rightarrow \P^-\C$ is
generally not surjective on the underlying sets but by
Theorem~\ref{calc_rules}, any element of $(\P^-\C)_{n-1}$ is a
composite of elements of $\P\C$ by only using the composition laws
$*_r$ for $r\geq 1$. Since the branching semi-globular nerve is
regular, then $\square_n^{gl^-}(x)$ is therefore equal to a sum of
elements of the form $\square_n^{gl^-}(h^-(z))$ where $z\in
(\P^-\C)_{n-1}$. Hence the surjectivity.  \epf

These three maps are likely to be injective for any non-contracting
$\omega$-category $\C$ but we have not yet been able to verify it.
The latter conjecture together with the thin elements conjecture means
that for $\omega$-categories modeling HDA, i.e. freely generated by
precubical sets, both simplicial homology theories $H_*^-$ and
$H_*^{gl^-}$ coincide.

Figure~\ref{recap} is a recapitulation of all constructions made so
far ($\mathcal{N}^{gl^+}$ is the merging semi-globular cut, that is
the Kan version of the merging nerve). The figure represents a
commutative diagram. The morphisms $h^-$ and $h^+$ are those defined
in \cite{sglob}.  The maps $\mathcal{N}^{-}\rightarrow
\mathcal{N}^{gl^-}$ and $\mathcal{N}^{+}\rightarrow
\mathcal{N}^{gl^+}$ are the maps constructed in Theorem~\ref{kani}.
The composites $\mathcal{N}^{gl}\rightarrow \mathcal{N}^{gl^-}$ and
$\mathcal{N}^{gl}\rightarrow \mathcal{N}^{gl^+}$ are induced by the
natural transformations  $\P\rightarrow \P^-$ (resp.
$\P\rightarrow \P^+$).

\section{Concluding discussion}

Any HDA can be modeled by a precubical set $K$. The prefix ``pre''
means that there are no degeneracy maps in the data.

We have worked so far with the strict globular $\omega$-category
$F(K)$ freely generated by the precubical set $K$. This paper
shows, hopefully with convincing arguments, that we can deal
directly with the non-contracting $\omega$-category
$\widehat{F(K)}$ obtained by making the path $\omega$-category
$\P F(K)$ a strict globular $\omega$-groupoid in a universal way,
without changing the information contained in the homology groups
introduced so far. The non-contracting $\omega$-category
$\widehat{F(K)}$ is obtained from the non-contracting
$\omega$-category ${F(K)}$ by adding inverses to morphisms of
dimension greater than $2$ and with respect to all composition
laws of dimension greater than $1$. This is very satisfactory
from the point of view of computer-scientific modeling because
there are no reasons for an homotopy between non-constant
execution paths to be not invertible.  In this new setting, the
globular nerve $\mathcal{N}^{gl}(\widehat{F(K)})$ of
$\widehat{F(K)}$ becomes a Kan complex.

Starting from a non-contracting $\omega$-category $\C$, we have then
introduced in this paper two $\omega$-categories $\P^-\C$ and $\P^+\C$
whose $0$-morphisms are the germs of $1$-morphisms of $\C$ beginning
(resp. ending) in the same way and whose higher dimensional morphisms
are the germs of homotopies between them. We have then obtained the
diagram of non-contracting $\omega$-categories
\[
\xymatrix{
& \P \C \ar@{->}[ld]_{h^-}\ar@{->}[rd]^{h^+}& \\
\P^-\C && \P^+\C
}
\]
Applied to the particular case $\C=\widehat{F(K)}$, one obtains
the diagram of strict globular $\omega$-groupoids
\[
\xymatrix{
& \P \widehat{F(K)} \ar@{->}[ld]_{h^-}\ar@{->}[rd]^{h^+}& \\
\P^-\widehat{F(K)} && \P^+\widehat{F(K)}
}
\]
whose corresponding simplicial homologies give us the globular
homology and the new branching and merging homologies.

Loosely speaking, the strict globular $\omega$-groupoid $\P
\widehat{F(K)}$ plays the role of the space of non-constant execution
paths of the HDA, and the strict globular $\omega$-groupoid $\P^-
\widehat{F(K)}$ (resp. $\P^+ \widehat{F(K)}$ ) plays the role of the
space of germs of non-constant execution paths beginning
(resp. ending) in the same way. So these results tell us that what
matters for the homological study of dihomotopy is to have a set of
states (the set of $0$-morphisms) and a space of non-constant
execution paths, the two other spaces being characterized by the
space of non-constant execution paths and the composition law $*_0$.
This idea is implemented  in \cite{flow} in a topological context.

\section{Acknowledgments}

I thank the referee for helpful comments and the editor for his
patience.

\textit{\\Institut de Recherche Math\'ematique Avanc\'ee\\
ULP et CNRS\\ 7 rue Ren\'e Descartes\\67084 Strasbourg Cedex\\France\\gaucher@irma.u-strasbg.fr}

\end{document}